\input amstex
\input amsppt.sty

\input texdraw

\magnification1200
\hsize14cm
\vsize19cm

\TagsOnRight

\def\BrGrAA{1}
\def\BresAO{2}
\def\FuKrAA{3}
\def\FuHaAA{4}
\def\KratBC{5}
\def\KupeAA{6}
\def\LascAZ{7}
\def\MacdAC{8}
\def\SagaAL{9}
\def\StanBI{10}

\def\AA{2.1}
\def\AB{2.2}
\def\ABa{2.4}
\def\ABb{2.6}
\def\AC{1.1}
\def\ACa{1.2}
\def\ACb{1.3}
\def\ACc{1.4}
\def\ACd{1.5}
\def\AD{5.1}
\def\AE{5.2}
\def\AF{5.3}
\def\AFa{5.5}
\def\AFb{5.4}
\def\AG{3.1}
\def\AGa{3.2}
\def\AGb{3.3}
\def\AH{4.1}
\def\AI{4.2}
\def\AJ{4.3}
\def\AK{4.4}
\def\AKa{2.3}
\def\AL{4.5}
\def\AMa{4.6}
\def\AM{4.7}
\def\AN{4.8}
\def\AO{4.9}
\def\AP{4.10}

\def\AR{2.7}
\def\AS{2.5}
\def\AT{6.1}
\def\AU{6.2}
\def\AV{6.3}
\def\AW{6.4}
\def\AX{6.5}
\def\AY{6.6}
\def\AZ{6.7}

\def\FA{1}
\def\FB{2}
\def\FC{3}
\def\FD{4}
\def\FE{5}


\def\RhombusA{\bsegment
  \rlvec(0.866025403784439 .5) \rlvec(0.866025403784439 -.5) 
  \rlvec(-0.866025403784439 -.5) \rlvec(-0.866025403784439 .5) 
  \savepos(0.866025403784439 -.5)(*ex *ey)
        \esegment
  \move(*ex *ey)
        }
\def\RhombusB{\bsegment
  \rlvec(0.866025403784439 .5) \rlvec(0 -1)
  \rlvec(-0.866025403784439 -.5) \rlvec(0 1) 
  \savepos(0 -1)(*ex *ey)
        \esegment
  \move(*ex *ey)
        }
\def\RhombusC{\bsegment
  \rlvec(0.866025403784439 -.5) \rlvec(0 -1)
  \rlvec(-0.866025403784439 .5) \rlvec(0 1) 
  \savepos(0.866025403784439 -.5)(*ex *ey)
        \esegment
  \move(*ex *ey)
        }

\def\hdSchritt{\bsegment
  \lpatt(.05 .13)
  \rlvec(0.866025403784439 -.5) 
  \savepos(0.866025403784439 -.5)(*ex *ey)
        \esegment
  \move(*ex *ey)
        }

\def\odaSchritt{\bsegment
  \lpatt(.05 .13)
  \rlvec(0.866025403784439 .5) 
  \savepos(0.866025403784439 .5)(*ex *ey)
        \esegment
  \move(*ex *ey)
        }

\def\ringerl(#1 #2){\move(#1 #2)\fcir f:0 r:.15}

\catcode`\@=11
\font\tenln    = line10
\font\tenlnw   = linew10

\newskip\Einheit \Einheit=0.5cm
\newcount\xcoord \newcount\ycoord
\newdimen\xdim \newdimen\ydim \newdimen\PfadD@cke \newdimen\Pfadd@cke

\newcount\@tempcnta
\newcount\@tempcntb

\newdimen\@tempdima
\newdimen\@tempdimb

\newdimen\@wholewidth
\newdimen\@halfwidth

\newcount\@xarg
\newcount\@yarg
\newcount\@yyarg
\newbox\@linechar
\newbox\@tempboxa
\newdimen\@linelen
\newdimen\@clnwd
\newdimen\@clnht

\newif\if@negarg

\def\@whilenoop#1{}
\def\@whiledim#1\do #2{\ifdim #1\relax#2\@iwhiledim{#1\relax#2}\fi}
\def\@iwhiledim#1{\ifdim #1\let\@nextwhile=\@iwhiledim
        \else\let\@nextwhile=\@whilenoop\fi\@nextwhile{#1}}

\def\@whileswnoop#1\fi{}
\def\@whilesw#1\fi#2{#1#2\@iwhilesw{#1#2}\fi\fi}
\def\@iwhilesw#1\fi{#1\let\@nextwhile=\@iwhilesw
         \else\let\@nextwhile=\@whileswnoop\fi\@nextwhile{#1}\fi}

\def\thinlines{\let\@linefnt\tenln \let\@circlefnt\tencirc
  \@wholewidth\fontdimen8\tenln \@halfwidth .5\@wholewidth}
\def\thicklines{\let\@linefnt\tenlnw \let\@circlefnt\tencircw
  \@wholewidth\fontdimen8\tenlnw \@halfwidth .5\@wholewidth}
\thinlines

\PfadD@cke1pt \Pfadd@cke0.5pt
\def\PfadDicke#1{\PfadD@cke#1 \divide\PfadD@cke by2 \Pfadd@cke\PfadD@cke \multiply\PfadD@cke by2}
\long\def\LOOP#1\REPEAT{\def\BODY{#1}\ITERATE}
\def\ITERATE{\BODY \let\next\ITERATE \else\let\next\relax\fi \next}
\let\REPEAT=\fi
\def\Punkt{\hbox{\raise-2pt\hbox to0pt{\hss$\ssize\bullet$\hss}}}
\def\DuennPunkt(#1,#2){\unskip
  \raise#2 \Einheit\hbox to0pt{\hskip#1 \Einheit
          \raise-2.5pt\hbox to0pt{\hss$\bullet$\hss}\hss}}
\def\NormalPunkt(#1,#2){\unskip
  \raise#2 \Einheit\hbox to0pt{\hskip#1 \Einheit
          \raise-3pt\hbox to0pt{\hss\twelvepoint$\bullet$\hss}\hss}}
\def\DickPunkt(#1,#2){\unskip
  \raise#2 \Einheit\hbox to0pt{\hskip#1 \Einheit
          \raise-4pt\hbox to0pt{\hss\fourteenpoint$\bullet$\hss}\hss}}
\def\Kreis(#1,#2){\unskip
  \raise#2 \Einheit\hbox to0pt{\hskip#1 \Einheit
          \raise-4pt\hbox to0pt{\hss\fourteenpoint$\circ$\hss}\hss}}

\def\Line@(#1,#2)#3{\@xarg #1\relax \@yarg #2\relax
\@linelen=#3\Einheit
\ifnum\@xarg =0 \@vline
  \else \ifnum\@yarg =0 \@hline \else \@sline\fi
\fi}

\def\@sline{\ifnum\@xarg< 0 \@negargtrue \@xarg -\@xarg \@yyarg -\@yarg
  \else \@negargfalse \@yyarg \@yarg \fi
\ifnum \@yyarg >0 \@tempcnta\@yyarg \else \@tempcnta -\@yyarg \fi
\ifnum\@tempcnta>6 \@badlinearg\@tempcnta0 \fi
\ifnum\@xarg>6 \@badlinearg\@xarg 1 \fi
\setbox\@linechar\hbox{\@linefnt\@getlinechar(\@xarg,\@yyarg)}%
\ifnum \@yarg >0 \let\@upordown\raise \@clnht\z@
   \else\let\@upordown\lower \@clnht \ht\@linechar\fi
\@clnwd=\wd\@linechar
\if@negarg \hskip -\wd\@linechar \def\@tempa{\hskip -2\wd\@linechar}\else
     \let\@tempa\relax \fi
\@whiledim \@clnwd <\@linelen \do
  {\@upordown\@clnht\copy\@linechar
   \@tempa
   \advance\@clnht \ht\@linechar
   \advance\@clnwd \wd\@linechar}%
\advance\@clnht -\ht\@linechar
\advance\@clnwd -\wd\@linechar
\@tempdima\@linelen\advance\@tempdima -\@clnwd
\@tempdimb\@tempdima\advance\@tempdimb -\wd\@linechar
\if@negarg \hskip -\@tempdimb \else \hskip \@tempdimb \fi
\multiply\@tempdima \@m
\@tempcnta \@tempdima \@tempdima \wd\@linechar \divide\@tempcnta \@tempdima
\@tempdima \ht\@linechar \multiply\@tempdima \@tempcnta
\divide\@tempdima \@m
\advance\@clnht \@tempdima
\ifdim \@linelen <\wd\@linechar
   \hskip \wd\@linechar
  \else\@upordown\@clnht\copy\@linechar\fi}

\def\@hline{\ifnum \@xarg <0 \hskip -\@linelen \fi
\vrule height\Pfadd@cke width \@linelen depth\Pfadd@cke
\ifnum \@xarg <0 \hskip -\@linelen \fi}

\def\@getlinechar(#1,#2){\@tempcnta#1\relax\multiply\@tempcnta 8
\advance\@tempcnta -9 \ifnum #2>0 \advance\@tempcnta #2\relax\else
\advance\@tempcnta -#2\relax\advance\@tempcnta 64 \fi
\char\@tempcnta}

\def\Vektor(#1,#2)#3(#4,#5){\unskip\leavevmode
  \xcoord#4\relax \ycoord#5\relax
      \raise\ycoord \Einheit\hbox to0pt{\hskip\xcoord \Einheit
         \Vector@(#1,#2){#3}\hss}}

\def\Vector@(#1,#2)#3{\@xarg #1\relax \@yarg #2\relax
\@tempcnta \ifnum\@xarg<0 -\@xarg\else\@xarg\fi
\ifnum\@tempcnta<5\relax
\@linelen=#3\Einheit
\ifnum\@xarg =0 \@vvector
  \else \ifnum\@yarg =0 \@hvector \else \@svector\fi
\fi
\else\@badlinearg\fi}

\def\@hvector{\@hline\hbox to 0pt{\@linefnt
\ifnum \@xarg <0 \@getlarrow(1,0)\hss\else
    \hss\@getrarrow(1,0)\fi}}

\def\@vvector{\ifnum \@yarg <0 \@downvector \else \@upvector \fi}

\def\@svector{\@sline
\@tempcnta\@yarg \ifnum\@tempcnta <0 \@tempcnta=-\@tempcnta\fi
\ifnum\@tempcnta <5
  \hskip -\wd\@linechar
  \@upordown\@clnht \hbox{\@linefnt  \if@negarg
  \@getlarrow(\@xarg,\@yyarg) \else \@getrarrow(\@xarg,\@yyarg) \fi}%
\else\@badlinearg\fi}

\def\@upline{\hbox to \z@{\hskip -.5\Pfadd@cke \vrule width \Pfadd@cke
   height \@linelen depth \z@\hss}}

\def\@downline{\hbox to \z@{\hskip -.5\Pfadd@cke \vrule width \Pfadd@cke
   height \z@ depth \@linelen \hss}}

\def\@upvector{\@upline\setbox\@tempboxa\hbox{\@linefnt\char'66}\raise
     \@linelen \hbox to\z@{\lower \ht\@tempboxa\box\@tempboxa\hss}}

\def\@downvector{\@downline\lower \@linelen
      \hbox to \z@{\@linefnt\char'77\hss}}

\def\@getlarrow(#1,#2){\ifnum #2 =\z@ \@tempcnta='33\else
\@tempcnta=#1\relax\multiply\@tempcnta \sixt@@n \advance\@tempcnta
-9 \@tempcntb=#2\relax\multiply\@tempcntb \tw@
\ifnum \@tempcntb >0 \advance\@tempcnta \@tempcntb\relax
\else\advance\@tempcnta -\@tempcntb\advance\@tempcnta 64
\fi\fi\char\@tempcnta}

\def\@getrarrow(#1,#2){\@tempcntb=#2\relax
\ifnum\@tempcntb < 0 \@tempcntb=-\@tempcntb\relax\fi
\ifcase \@tempcntb\relax \@tempcnta='55 \or
\ifnum #1<3 \@tempcnta=#1\relax\multiply\@tempcnta
24 \advance\@tempcnta -6 \else \ifnum #1=3 \@tempcnta=49
\else\@tempcnta=58 \fi\fi\or
\ifnum #1<3 \@tempcnta=#1\relax\multiply\@tempcnta
24 \advance\@tempcnta -3 \else \@tempcnta=51\fi\or
\@tempcnta=#1\relax\multiply\@tempcnta
\sixt@@n \advance\@tempcnta -\tw@ \else
\@tempcnta=#1\relax\multiply\@tempcnta
\sixt@@n \advance\@tempcnta 7 \fi\ifnum #2<0 \advance\@tempcnta 64 \fi
\char\@tempcnta}

\def\Diagonale(#1,#2)#3{\unskip\leavevmode
  \xcoord#1\relax \ycoord#2\relax
      \raise\ycoord \Einheit\hbox to0pt{\hskip\xcoord \Einheit
         \Line@(1,1){#3}\hss}}
\def\AntiDiagonale(#1,#2)#3{\unskip\leavevmode
  \xcoord#1\relax \ycoord#2\relax 
      \raise\ycoord \Einheit\hbox to0pt{\hskip\xcoord \Einheit
         \Line@(1,-1){#3}\hss}}
\def\Pfad(#1,#2),#3\endPfad{\unskip\leavevmode
  \xcoord#1 \ycoord#2 \thicklines\ZeichnePfad#3\endPfad\thinlines}
\def\ZeichnePfad#1{\ifx#1\endPfad\let\next\relax
  \else\let\next\ZeichnePfad
    \ifnum#1=1
      \raise\ycoord \Einheit\hbox to0pt{\hskip\xcoord \Einheit
         \vrule height\Pfadd@cke width1 \Einheit depth\Pfadd@cke\hss}%
      \advance\xcoord by 1
    \else\ifnum#1=2
      \raise\ycoord \Einheit\hbox to0pt{\hskip\xcoord \Einheit
        \hbox{\hskip-\PfadD@cke\vrule height1 \Einheit width\PfadD@cke depth0pt}\hss}%
      \advance\ycoord by 1
    \else\ifnum#1=3
      \raise\ycoord \Einheit\hbox to0pt{\hskip\xcoord \Einheit
         \Line@(1,1){1}\hss}
      \advance\xcoord by 1
      \advance\ycoord by 1
    \else\ifnum#1=4
      \raise\ycoord \Einheit\hbox to0pt{\hskip\xcoord \Einheit
         \Line@(1,-1){1}\hss}
      \advance\xcoord by 1
      \advance\ycoord by -1
    \fi\fi\fi\fi
  \fi\next}
\def\hSSchritt{\leavevmode\raise-.4pt\hbox to0pt{\hss.\hss}\hskip.2\Einheit
  \raise-.4pt\hbox to0pt{\hss.\hss}\hskip.2\Einheit
  \raise-.4pt\hbox to0pt{\hss.\hss}\hskip.2\Einheit
  \raise-.4pt\hbox to0pt{\hss.\hss}\hskip.2\Einheit
  \raise-.4pt\hbox to0pt{\hss.\hss}\hskip.2\Einheit}
\def\vSSchritt{\vbox{\baselineskip.2\Einheit\lineskiplimit0pt
\hbox{.}\hbox{.}\hbox{.}\hbox{.}\hbox{.}}}
\def\DSSchritt{\leavevmode\raise-.4pt\hbox to0pt{%
  \hbox to0pt{\hss.\hss}\hskip.2\Einheit
  \raise.2\Einheit\hbox to0pt{\hss.\hss}\hskip.2\Einheit
  \raise.4\Einheit\hbox to0pt{\hss.\hss}\hskip.2\Einheit
  \raise.6\Einheit\hbox to0pt{\hss.\hss}\hskip.2\Einheit
  \raise.8\Einheit\hbox to0pt{\hss.\hss}\hss}}
\def\dSSchritt{\leavevmode\raise-.4pt\hbox to0pt{%
  \hbox to0pt{\hss.\hss}\hskip.2\Einheit
  \raise-.2\Einheit\hbox to0pt{\hss.\hss}\hskip.2\Einheit
  \raise-.4\Einheit\hbox to0pt{\hss.\hss}\hskip.2\Einheit
  \raise-.6\Einheit\hbox to0pt{\hss.\hss}\hskip.2\Einheit
  \raise-.8\Einheit\hbox to0pt{\hss.\hss}\hss}}
\def\SPfad(#1,#2),#3\endSPfad{\unskip\leavevmode
  \xcoord#1 \ycoord#2 \ZeichneSPfad#3\endSPfad}
\def\ZeichneSPfad#1{\ifx#1\endSPfad\let\next\relax
  \else\let\next\ZeichneSPfad
    \ifnum#1=1
      \raise\ycoord \Einheit\hbox to0pt{\hskip\xcoord \Einheit
         \hSSchritt\hss}%
      \advance\xcoord by 1
    \else\ifnum#1=2
      \raise\ycoord \Einheit\hbox to0pt{\hskip\xcoord \Einheit
        \hbox{\hskip-2pt \vSSchritt}\hss}%
      \advance\ycoord by 1
    \else\ifnum#1=3
      \raise\ycoord \Einheit\hbox to0pt{\hskip\xcoord \Einheit
         \DSSchritt\hss}
      \advance\xcoord by 1
      \advance\ycoord by 1
    \else\ifnum#1=4
      \raise\ycoord \Einheit\hbox to0pt{\hskip\xcoord \Einheit
         \dSSchritt\hss}
      \advance\xcoord by 1
      \advance\ycoord by -1
    \fi\fi\fi\fi
  \fi\next}
\def\Koordinatenachsen(#1,#2){\unskip
 \hbox to0pt{\hskip-.5pt\vrule height#2 \Einheit width.5pt depth1 \Einheit}%
 \hbox to0pt{\hskip-1 \Einheit \xcoord#1 \advance\xcoord by1
    \vrule height0.25pt width\xcoord \Einheit depth0.25pt\hss}}
\def\Koordinatenachsen(#1,#2)(#3,#4){\unskip
 \hbox to0pt{\hskip-.5pt \ycoord-#4 \advance\ycoord by1
    \vrule height#2 \Einheit width.5pt depth\ycoord \Einheit}%
 \hbox to0pt{\hskip-1 \Einheit \hskip#3\Einheit
    \xcoord#1 \advance\xcoord by1 \advance\xcoord by-#3
    \vrule height0.25pt width\xcoord \Einheit depth0.25pt\hss}}
\def\Gitter(#1,#2){\unskip \xcoord0 \ycoord0 \leavevmode
  \LOOP\ifnum\ycoord<#2
    \loop\ifnum\xcoord<#1
      \raise\ycoord \Einheit\hbox to0pt{\hskip\xcoord \Einheit\Punkt\hss}%
      \advance\xcoord by1
    \repeat
    \xcoord0
    \advance\ycoord by1
  \REPEAT}
\def\Gitter(#1,#2)(#3,#4){\unskip \xcoord#3 \ycoord#4 \leavevmode
  \LOOP\ifnum\ycoord<#2
    \loop\ifnum\xcoord<#1
      \raise\ycoord \Einheit\hbox to0pt{\hskip\xcoord \Einheit\Punkt\hss}%
      \advance\xcoord by1
    \repeat
    \xcoord#3
    \advance\ycoord by1
  \REPEAT}
\def\Label#1#2(#3,#4){\unskip \xdim#3 \Einheit \ydim#4 \Einheit
  \def\lo{\advance\xdim by-.5 \Einheit \advance\ydim by.5 \Einheit}%
  \def\llo{\advance\xdim by-.25cm \advance\ydim by.5 \Einheit}%
  \def\loo{\advance\xdim by-.5 \Einheit \advance\ydim by.25cm}%
  \def\o{\advance\ydim by.25cm}%
  \def\ro{\advance\xdim by.5 \Einheit \advance\ydim by.5 \Einheit}%
  \def\rro{\advance\xdim by.25cm \advance\ydim by.5 \Einheit}%
  \def\roo{\advance\xdim by.5 \Einheit \advance\ydim by.25cm}%
  \def\l{\advance\xdim by-.30cm}%
  \def\r{\advance\xdim by.30cm}%
  \def\lu{\advance\xdim by-.5 \Einheit \advance\ydim by-.6 \Einheit}%
  \def\llu{\advance\xdim by-.25cm \advance\ydim by-.6 \Einheit}%
  \def\luu{\advance\xdim by-.5 \Einheit \advance\ydim by-.30cm}%
  \def\u{\advance\ydim by-.30cm}%
  \def\ru{\advance\xdim by.5 \Einheit \advance\ydim by-.6 \Einheit}%
  \def\rru{\advance\xdim by.25cm \advance\ydim by-.6 \Einheit}%
  \def\ruu{\advance\xdim by.5 \Einheit \advance\ydim by-.30cm}%
  #1\raise\ydim\hbox to0pt{\hskip\xdim
     \vbox to0pt{\vss\hbox to0pt{\hss$#2$\hss}\vss}\hss}%
}

\font@\twelverm=cmr10 scaled\magstep1
\font@\twelveit=cmti10 scaled\magstep1
\font@\twelvebf=cmbx10 scaled\magstep1
\font@\twelvei=cmmi10 scaled\magstep1
\font@\twelvesy=cmsy10 scaled\magstep1
\font@\twelveex=cmex10 scaled\magstep1

\newtoks\twelvepoint@
\def\twelvepoint{\normalbaselineskip15\p@
 \abovedisplayskip15\p@ plus3.6\p@ minus10.8\p@
 \belowdisplayskip\abovedisplayskip
 \abovedisplayshortskip\z@ plus3.6\p@
 \belowdisplayshortskip8.4\p@ plus3.6\p@ minus4.8\p@
 \textonlyfont@\rm\twelverm \textonlyfont@\it\twelveit
 \textonlyfont@\sl\twelvesl \textonlyfont@\bf\twelvebf
 \textonlyfont@\smc\twelvesmc \textonlyfont@\tt\twelvett
%
 \ifsyntax@ \def\big##1{{\hbox{$\left##1\right.$}}}%
  \let\Big\big \let\bigg\big \let\Bigg\big
 \else
  \textfont\z@=\twelverm  \scriptfont\z@=\tenrm  \scriptscriptfont\z@=\sevenrm
  \textfont\@ne=\twelvei  \scriptfont\@ne=\teni  \scriptscriptfont\@ne=\seveni
  \textfont\tw@=\twelvesy \scriptfont\tw@=\tensy \scriptscriptfont\tw@=\sevensy
  \textfont\thr@@=\twelveex \scriptfont\thr@@=\tenex
        \scriptscriptfont\thr@@=\tenex
  \textfont\itfam=\twelveit \scriptfont\itfam=\tenit
        \scriptscriptfont\itfam=\tenit
  \textfont\bffam=\twelvebf \scriptfont\bffam=\tenbf
        \scriptscriptfont\bffam=\sevenbf
  \setbox\strutbox\hbox{\vrule height10.2\p@ depth4.2\p@ width\z@}%
  \setbox\strutbox@\hbox{\lower.6\normallineskiplimit\vbox{%
        \kern-\normallineskiplimit\copy\strutbox}}%
 \setbox\z@\vbox{\hbox{$($}\kern\z@}\bigsize@=1.4\ht\z@
 \fi
 \normalbaselines\rm\ex@.2326ex\jot3.6\ex@\the\twelvepoint@}

\font@\fourteenrm=cmr10 scaled\magstep2
\font@\fourteenit=cmti10 scaled\magstep2
\font@\fourteensl=cmsl10 scaled\magstep2
\font@\fourteensmc=cmcsc10 scaled\magstep2
\font@\fourteentt=cmtt10 scaled\magstep2
\font@\fourteenbf=cmbx10 scaled\magstep2
\font@\fourteeni=cmmi10 scaled\magstep2
\font@\fourteensy=cmsy10 scaled\magstep2
\font@\fourteenex=cmex10 scaled\magstep2
\font@\fourteenmsa=msam10 scaled\magstep2
\font@\fourteeneufm=eufm10 scaled\magstep2
\font@\fourteenmsb=msbm10 scaled\magstep2
\newtoks\fourteenpoint@
\def\fourteenpoint{\normalbaselineskip15\p@
 \abovedisplayskip18\p@ plus4.3\p@ minus12.9\p@
 \belowdisplayskip\abovedisplayskip
 \abovedisplayshortskip\z@ plus4.3\p@
 \belowdisplayshortskip10.1\p@ plus4.3\p@ minus5.8\p@
 \textonlyfont@\rm\fourteenrm \textonlyfont@\it\fourteenit
 \textonlyfont@\sl\fourteensl \textonlyfont@\bf\fourteenbf
 \textonlyfont@\smc\fourteensmc \textonlyfont@\tt\fourteentt
%
 \ifsyntax@ \def\big##1{{\hbox{$\left##1\right.$}}}%
  \let\Big\big \let\bigg\big \let\Bigg\big
 \else
  \textfont\z@=\fourteenrm  \scriptfont\z@=\twelverm  \scriptscriptfont\z@=\tenrm
  \textfont\@ne=\fourteeni  \scriptfont\@ne=\twelvei  \scriptscriptfont\@ne=\teni
  \textfont\tw@=\fourteensy \scriptfont\tw@=\twelvesy \scriptscriptfont\tw@=\tensy
  \textfont\thr@@=\fourteenex \scriptfont\thr@@=\twelveex
        \scriptscriptfont\thr@@=\twelveex
  \textfont\itfam=\fourteenit \scriptfont\itfam=\twelveit
        \scriptscriptfont\itfam=\twelveit
  \textfont\bffam=\fourteenbf \scriptfont\bffam=\twelvebf
        \scriptscriptfont\bffam=\tenbf
  \setbox\strutbox\hbox{\vrule height12.2\p@ depth5\p@ width\z@}%
  \setbox\strutbox@\hbox{\lower.72\normallineskiplimit\vbox{%
        \kern-\normallineskiplimit\copy\strutbox}}%
 \setbox\z@\vbox{\hbox{$($}\kern\z@}\bigsize@=1.7\ht\z@
 \fi
 \normalbaselines\rm\ex@.2326ex\jot4.3\ex@\the\fourteenpoint@}

\catcode`\@=13

\def\la{\lambda}
\def\so{\operatorname{\text {\it so}}}

\def\sp{\operatorname{\text {\it sp}}}
\def\PP{\operatorname{\text {\it PP}}}
\def\TCPP{\operatorname{\text {\it TCPP}}}
\def\SPP{\operatorname{\text {\it SPP}}}
\def\C{\Bbb C}
\def\x{\bold x}
\def\oddcols{\operatorname{oddcols}}

\topmatter
\title A factorization theorem for classical group characters, with
applications to plane partitions and rhombus tilings
\endtitle
\author M.~Ciucu and C.~Krattenthaler
\endauthor
\affil Department of Mathematics,
Indiana University,\\
Bloomington, IN 47405-5701, USA\\\vskip6pt
Fakult\"at f\"ur Mathematik der Universit\"at Wien,\\
Nordbergstra{\ss}e 15, A-1090 Wien, Austria.\\
WWW: \tt http://www.mat.univie.ac.at/\~{}kratt
\endaffil
\address Department of Mathematics, Indiana University, Bloomington, IN 47405-5701,
USA
\endaddress
\address Fakult\"at f\"ur Mathematik der Universit\"at Wien,
Nordbergstra\ss e 15, A-1090 Wien, Austria.
\endaddress
\thanks{The research of the first author was
partially supported by NSF grant DMS-0500616. The
research of the second author was partially supported by the Austrian
Science Foundation FWF, grant S9607-N13,
in the framework of the National Research Network
``Analytic Combinatorics and Probabilistic Number Theory."
This work was done during the authors' stay at the
Erwin Schr\"odinger Institute for Physics and Mathematics,
Vienna, during the programme ``Combinatorics and Statistical Physics''
in Spring~2008.}\endthanks
\abstract
We prove that a Schur function of rectangular shape $(M^n)$ whose
variables are specialized to $x_1,x_1^{-1},\dots,x_n,x_n^{-1}$
factorizes into a product of two odd orthogonal characters of
rectangular shape, one of which is evaluated at $-x_1,\dots,-x_n$, if
$M$ is even, while it factorizes into a product of a symplectic
character and an even orthogonal
character, both of rectangular shape, if $M$ is odd. It is furthermore shown
that the first factorization implies a factorization theorem for
rhombus tilings of a hexagon, which has an equivalent formulation in
terms of plane partitions. A similar factorization theorem is proven
for the sum of two Schur functions of respective rectangular shapes
$(M^n)$ and $(M^{n-1})$.
\endabstract
\endtopmatter
\document

\leftheadtext{M. Ciucu and C. Krattenthaler}
\rightheadtext{A factorization theorem for classical group characters}

\subhead 1. Introduction\endsubhead
The purpose of this note is to prove curious factorization properties
for Schur functions of rectangular shape, which seem to have escaped
the attention of previous authors. (We refer the reader to Section~2
for all definitions.) More precisely, we show that a Schur function of
rectangular shape $(M^n)$ which is evaluated at $x_1,x_2,\dots,x_n$
and their reciprocals $x_1^{-1},x_2^{-1},\dots,x_n^{-1}$ factorizes
into two factors, and the same is true for the sum of two Schur functions
of respective shapes $(M^n)$ and $(M^{n-1})$. 

We begin by describing explicitly the case of one Schur function.
If $M$ is even, then both factors are odd
orthogonal characters of rectangular shape, one of them evaluated at the
variables $x_1,x_2,\dots,x_n$, but the other is evaluated
at $-x_1,-x_2,\dots,-x_n$. If $M$ is odd, then one factor
is a symplectic character of rectangular shape, while the other is an
even orthogonal character of rectangular shape, both being evaluated at
$x_1,x_2,\dots,x_n$. 
The case of even $M$ of this factorization
property is presented in the following theorem.

\proclaim{Theorem 1}For any non-negative integers $m$ and $n$, we have
$$\multline 
s_{((2m)^n)}(x_1,x_1^{-1},x_2,x_2^{-1},\dots,x_n,x_n^{-1})\\
=(-1)^{mn}
\so_{(m^n)}(x_1,x_2,\dots,x_n)\,
\so_{(m^n)}(-x_1,-x_2,\dots,-x_n).
\endmultline
\tag\AC$$
\endproclaim

If $M$ is odd, then the factorization takes the following form.

\proclaim{Theorem 2}For any non-negative integers $m$ and $n$, we have
$$\multline 
s_{((2m+1)^n)}(x_1,x_1^{-1},x_2,x_2^{-1},\dots,x_n,x_n^{-1})\\
=
\sp_{(m^n)}(x_1,x_2,\dots,x_n)\,
o^{even}_{((m+1)^n)}(x_1,x_2,\dots,x_n).
\endmultline
\tag\ACa$$
\endproclaim

Since our identities involve classical group characters, one might
ask whether there are representation-theoretic interpretations of
these identities. At first sight, this seems to be a
difficult question because of the somewhat ``incoherent" 
right-hand sides of (\AC) and (\ACa). However, there is a uniform way of
writing the factorization identities of Theorems~1 and 2 that was
pointed out to us by Soichi Okada. Namely, 
by comparing (\AB) and (\ABb), and by using (\AKa)
and (\AR), we see that
$$(-1)^{\sum _{i=1} ^{N}\la_i}\so_\la(-x_1,-x_2,\dots,-x_N)
\prod _{i=1} ^{N}(x_i^{1/2}+x_i^{-1/2})=
o^{even}_{\la+\frac {1} {2}}(x_1,x_2\dots,x_N),
\tag\ACb
$$
where $\la+\frac {1} {2}$ is short for $(\la_1+\frac {1} {2},\la_2+\frac {1}
{2},\dots,\la_N+\frac {1} {2})$. Furthermore, by 
comparing (\AB) and (\ABa), and by using (\AKa)
and (\AS), we see that
$$\sp_\la(x_1,x_2,\dots,x_N)
\prod _{i=1} ^{N}(x_i^{1/2}+x_i^{-1/2})=
\so_{\la+\frac {1} {2}}(x_1,x_2\dots,x_N).
\tag\ACc
$$
Theorem~1 and 2 may therefore be uniformly stated as
$$\multline 
\bigg(\prod _{i=1} ^{n}(x_i^{1/2}+x_i^{-1/2})\bigg)
s_{(M^n)}(x_1,x_1^{-1},x_2,x_2^{-1},\dots,x_n,x_n^{-1})\\
=\so_{\big((\frac {M} {2})^n\big)}(x_1,x_2,\dots,x_n)\,
o^{even}_{\big((\frac {M+1} {2})^n\big)}(x_1,x_2,\dots,x_n).
\endmultline
\tag\ACd$$
Written in this way, this may, in the end, lead to
a representation-theoretic interpretation of these identities,
although we must confess that we are not able to offer such an 
interpretation.

On the other hand, we are able to offer a {\it combinatorial\/} 
interpretation for Theorem~1. 
As we show in Section~5, if we specialize $x_1=x_2=\dots=x_n=1$ in
Theorem~1, then one obtains a factorization theorem for rhombus
tilings of a hexagon, which has also a natural, equivalent
formulation as a factorization theorem for plane partitions (see (\AD)). 
It is, in fact, this factorization theorem which we observed first,
and which formed the starting point of this
work. We suspect that a more general factorization theorem for rhombus
tilings is lurking behind. If there is a natural combinatorial 
interpretation of Theorem~2 is less clear. We make an attempt, also
in Section~5, but we consider it not entirely satisfactory.

The case of the sum of two Schur functions of rectangular shapes is
treated in Section~6. Namely, we show that there are very similar
factorization theorems for the sum of two Schur functions of 
respective rectangular shapes $(M^n)$ and $(M^{n-1})$
(see Theorems~3 and 4). The existence
of these was pointed out to us by Ron King.

The proofs of Theorems~1 and 2 are given in Section~4. 
These proofs are based on an auxiliary identity which is established in 
Section~3 (see Lemma~1). The proofs of Theorems~3 and 4 are based on
two further auxiliary identities of similar type, which are also
presented and proved in Section~3 (see Lemmas~2 and 3).

\subhead 2. Classical group characters\endsubhead
In this section we recall the definitions of the classical group
characters involved in the factorizations in Theorems~1 and 2.
We also briefly touch upon their significance in representation
theory.

Given a partition $\la=(\la_1,\la_2,\dots,\la_N)$ 
(i.e., a non-increasing sequence of non-negative integers) 
the {\it Schur function} $s_\la(x_1,x_2,\dots,x_N)$ is defined
by (see \cite{\FuHaAA, p.~403, (A.4)}, \cite{\LascAZ, Prop.~1.4.4},
or \cite{\MacdAC, Ch.~I, (3.1)})
$$
s_\lambda (x_1, x_2,\dots ,x_N)=\frac {
 \det\limits _{1\le h,t\le N}(x_h^{\lambda _t+N-t})} {
\det\limits_{1\le h,t\le N}(x_h^{N-t})}.
\tag\AA
$$
It is not difficult to see that the denominator in (\AA) 
cancels out, so that any Schur function
$s_\la(x_1,x_2,\dots,x_N)$ is in fact a {\it polynomial\/} in
$x_1,x_2,\dots,x_N$, and is thus well-defined for any choice of the
variables $x_1,x_2,\dots,x_N$. It is well-known (cf\. 
\cite{\FuHaAA, \S24.2}) that
$s_\la(x_1,x_2,\dots,x_N)$ is an irreducible
character of $SL_N(\C)$ (respectively $GL_N(\C)$).

Given a non-increasing sequence
$\la=(\la_1,\la_2,\dots,\la_n)$ of integers or
half-integers (the latter being, by definition, positive odd integers
divided by 2), the {\it odd orthogonal character} $\so_\la(x_1,x_2,\dots,x_N)$
is defined by (see \cite{\FuHaAA, (24.28)})
$$
\so_\la(x_1,x_2,\dots,x_N)=\frac {\det\limits_{1\le h,t\le
N}(x_h^{\la_t+N-t+\frac 12}-x_h^{-(\la_t+N-t+\frac 12)})}
{\det\limits_{1\le h,t\le
N}(x_h^{N-t+\frac 12}-x_h^{-(N-t+\frac 12)})}.
\tag\AB
$$
Again, it is not difficult to see that the denominator in (\AB) 
cancels out, so that any odd orthogonal character 
$\so_\la(x_1,x_2,\dots,x_N)$ is in fact a {\it Laurent polynomial\/} in
$x_1,x_2,\dots,x_N$ (i.e., a polynomial in 
$x_1,x_1^{-1},x_2,x_2^{-1},\dots,x_N,x_N^{-1}$),
and is thus well-defined for any choice of the
variables $x_1,x_2,\dots,x_N$ such that all of them are non-zero. 
By the {\it Weyl denominator formula for
type $B$} (cf.~\cite{\FuHaAA, Lemma~24.3, Ex.~A.62}),
$$\multline 
\det\limits_{1\le h,t\le
n}\left(x_h^{n-t+\frac 12}-x_h^{-(n-t-\frac 12)}\right)\\=
(x_1x_2\cdots x_n)^{-n+\frac 12}
\prod _{1\le h<t\le n} ^{}(x_h-x_t)(x_hx_t-1)\prod _{h=1}
^{n}(x_h-1),
\endmultline\tag\AKa$$
the denominator in (\AB) can be actually evaluated in product form.
It is well-known (cf\. \cite{\FuHaAA, \S24.2}) that
$so_\la(x_1,x_2,\dots,x_N)$ is an irreducible
character of $SO_{2N+1}(\C)$.

Given a partition 
$\la=(\la_1,\la_2,\dots,\la_n)$, 
the {\it symplectic character} $\sp_\la(x_1,x_2,\dots,x_N)$ is
defined by (see \cite{\FuHaAA, (24.18)})
$$
\sp_\la(x_1,x_2,\dots,x_N)
=\frac {\det\limits_{1\le h,t\le
N}(x_h^{\la_t+N-t+1}-x_h^{-(\la_t+N-t+1)})} 
{\det\limits_{1\le h,t\le
N}(x_h^{N-t+1}-x_h^{-(N-t+1)})}.
\tag\ABa$$
Similarly to odd orthogonal characters,
$\sp_\la(x_1,x_2,\dots,x_N)$ is a {\it Laurent polynomial\/} in
$x_1,x_2,\mathbreak\dots,x_N$,
and is thus well-defined for any choice of the
variables $x_1,x_2,\dots,x_N$ such that all of them are non-zero. 
By the {\it Weyl denominator formula for type~$C$} 
(cf.~\cite{\FuHaAA, Lemma~24.3, Ex.~A.52}),
$$
\det_{1\le h,t\le n}(x_h^{n-t+1}-x_h^{n-t+1})=(x_1\cdots x_n)^{-n}\prod
_{1\le h<t\le n} ^{}(x_h-x_t)(x_hx_t-1)
\prod _{h=1} ^{n}(x_h^2-1),
\tag\AS$$
the denominator in (\ABa) can be actually evaluated in product form.
Furthermore,\linebreak 
$sp_\la(x_1,x_2,\dots,x_N)$ is an irreducible
character of $Sp_{2N}(\C)$ (cf\. \cite{\FuHaAA, \S24.2}).

Finally, given a non-increasing sequence 
$\la=(\la_1,\la_2,\dots,\la_N)$ of {\it positive} integers or half-integers, 
the {\it even orthogonal character} $o^{even}_\la(x_1,x_2,\dots,x_N)$ is
given by (see \cite{\FuHaAA, (24.40) plus the remarks on top of page~411})
$$
o^{even}_\la(x_1,x_2,\dots,x_N)
=2\frac {\det\limits_{1\le h,t\le
N}(x_h^{\la_t+N-t}+x_h^{-(\la_t+N-t)})} 
{\det\limits_{1\le h,t\le
N}(x_h^{N-t}+x_h^{-(N-t)})}.
\tag\ABb$$
Here as well,
$o^{even}_\la(x_1,x_2,\dots,x_N)$ is a {\it Laurent polynomial\/} in
$x_1,x_2,\dots,x_N$,
and is thus well-defined for any choice of the
variables $x_1,x_2,\dots,x_N$ such that all of them are non-zero.
the {\it Weyl denominator formula for type~$D$} 
(cf.~\cite{\FuHaAA, Lemma~24.3, Ex.~A.66}),
$$
\det_{1\le h,t\le n}(x_h^{n-t}+x_h^{-(n-t)})=
2\cdot(x_1\cdots x_n)^{-n+1}\prod
_{1\le h<t\le n} ^{}(x_h-x_t)(x_hx_t-1),
\tag\AR$$
the denominator in (\ABb) can be again evaluated in product form.
It is well-known (cf\. \cite{\FuHaAA, \S24.2}) that
$o^{even}_\la(x_1,x_2,\dots,x_N)$ is an irreducible
character of $O_{2N}(\C)$. When restricted to $SO_{2N}(\C)$, it
splits into two different irreducible characters of $SO_{2N}(\C)$.
(The reader should observe that we assumed $\la_N>0$. If we had
allowed $\la_N=0$, then we would have to divide the right-hand
side of (\ABb) by 2 in order to obtain an irreducible character of
$O_{2N}(\C)$ or $SO_{2N}(\C)$.)

\subhead 3. Auxiliary identities\endsubhead
The proofs of Theorems~1 and 2, and those of Theorems~3 and 4 in
Section~6, hinge upon certain multivariable
identities which we prove in this section. Before we are able to state
these identities, we have to introduce some notation.

Throughout, we use the standard notation $[n]:=\{1,2,\dots,n\}$.
Let $A$ and $B$ be given subsets of the set of positive integers.
Slightly abusing notation for resultants from \cite{\LascAZ}, we define
$$R\left(A,B^{-1}\right):=\prod _{a\in A} ^{}\prod _{b\in B} ^{}(x_a-x_b^{-1}).$$
(In order to avoid any confusion on the part of the reader: the
symbol $B^{-1}$ in $R(A,B^{-1})$ has no meaning by itself, the
exponent $-1$ is just there to indicate that the reciprocals of the
variables indexed by $B$ are used on the right-hand side of the
definition.)
Furthermore, we define
$$V(A):=\underset{a<b}\to{\prod _{a,b\in A} ^{}}(x_a-x_b)$$
and
$$V\!\left(A^{-1}\right):=\underset{a<b}\to{\prod _{a,b\in A} ^{}}(x_a^{-1}-x_b^{-1}).$$
In all these definitions, empty products have to be interpreted as $1$.

Now we are in the position to state and prove the identity on which
the proofs of Theorems~1 and 2 are based.

\proclaim{Lemma 1}
For all positive integers $N$, there holds the identity
$$\multline 
\underset{\vert A\vert=N}\to{\sum _{A\subseteq[2N]} ^{}}
V(A)V\!\left(A^{-1}\right)
V(A^c)V\!\left((A^c)^{-1}\right)R\!\left(A,A^{-1}\right)
R\!\left(A^c,(A^c)^{-1}\right)
\\=
\sum _{A\subseteq[2N]} ^{}
V(A)V\!\left(A^{-1}\right)
V(A^c)V\!\left((A^c)^{-1}\right)R\!\left(A,(A^c)^{-1}\right)
R\!\left(A^c,A^{-1}\right),
\endmultline\tag\AG$$
where $A^c$ denotes the complement of $A$ in $[2N]$.
\endproclaim
\demo{Proof} 
We prove the assertion by induction on $N$. For $N=1$, identity~(\AG)
reduces to
$$
2(x_1-x_1^{-1})(x_2-x_2^{-1})=
2(x_1-x_2)(x_1^{-1}-x_2^{-1})+
2(x_1-x_2^{-1})(x_2-x_1^{-1}),
$$
which can be readily verified.

Now let us suppose that we have proved (\AG) with $N$ replaced by
$N-1$. Below, we shall show that this implies that the identity (\AG) is
true if we specialize $x_1$ to $x_2$. Let us for the moment
suppose that this is already done. Since, as is easy to
see, both sides of (\AG) are symmetric in the variables
$x_1,x_2,\dots,x_{2N}$, as well as they remain invariant up to an overall
multiplicative sign if we replace $x_1$ by $x_1^{-1}$
(on the left-hand side of (\AG), the latter assertion
already applies to each individual summand; 
on the right-hand side, one has to group the summands in pairs:
the summands corresponding to a set $A$ and the symmetric difference 
$A\triangle \{1\}$ have to be considered together), this means
that identity (\AG) holds for {\it all\/} specializations of $x_1$ to
one of $x_2,x_2^{-1},x_3,x_3^{-1},\dots,x_{2N},x_{2N}^{-1}$. 
These are $4N-2$
specializations. On the other hand, as Laurent polynomials in $x_1$, 
the degree of both sides of (\AG) is at most $2N-1$. (That is, the
maximal exponent $e$ of a power $x_1^e$ is $e=2N-1$, and the minimal
exponent is $e=-(2N-1)$.) Hence, both sides of (\AG) must agree up to
a multiplicative constant. In order to determine this multiplicative
constant, we compare coefficients of
$$x_1^{2N-1}x_2^{2N-3}\cdots x_N^1\cdot
x_{N+1}^{2N-1}x_{N+2}^{2N-3}\cdots x_{2N}^1
$$
on both sides of (\AG). On the left-hand side, the only terms
contributing are the ones corresponding to $A=\{1,2,\dots,N\}$ and
$A=\{N+1,N+2,\dots,2N\}$ (in fact, the corresponding terms are
equal to each other), both of them contributing a coefficient of 
$1$. On the right-hand side, the only terms
contributing are again the ones corresponding to $A=\{1,2,\dots,N\}$ and
$A=\{N+1,N+2,\dots,2N\}$ (being
equal to each other), both of them also contributing a coefficient of 
$1$.

It remains to prove that, under the induction hypothesis,
Equation~(\AG) holds for $x_1=x_2$. Indeed, under
this specialization, terms corresponding to sets $A$ which contain
{\it both\/} $1$ and $2$ and to those which contain {\it neither}  
$1$ {\it nor} $2$ vanish on both sides of (\AG), because of the
appearance of the Vandermonde products $V(A)$ respectively $V(A^c)$.
Therefore, if $x_1=x_2$, identity (\AG) reduces to
$$\multline 
(x_1-x_1^{-1})^2\bigg(\prod _{j=3}
^{2N}(x_1-x_j)(x_1^{-1}-x_j^{-1})(x_1-x_j^{-1})(x_j-x_1^{-1})\bigg)\\
\times
\underset{\vert A\vert=N-1}\to{\sum _{A\subseteq \{3,4,\dots,2N\}} ^{}}
V(A)V\!\left(A^{-1}\right)
V(A^c)V\!\left((A^c)^{-1}\right)R\!\left(A,A^{-1}\right)
R\!\left(A^c,(A^c)^{-1}\right)
\\=
(x_1-x_1^{-1})^2\bigg(\prod _{j=3}
^{2N}(x_1-x_j)(x_1^{-1}-x_j^{-1})(x_1-x_j^{-1})(x_j-x_1^{-1})\bigg)\\
\times
\sum _{A\subseteq \{3,4,\dots,2N\}} ^{}
V(A)V\!\left(A^{-1}\right)
V(A^c)V\!\left((A^c)^{-1}\right)R\!\left(A,(A^c)^{-1}\right)
R\!\left(A^c,A^{-1}\right).
\endmultline$$
After clearing the product which is common to both sides, we see that
the remaining identity is equivalent to (\AG) with $N$ replaced by
$N-1$, the latter being true due to the induction hypothesis. This
finishes the proof of the lemma.\quad \quad \qed
\enddemo

On the other hand, the proofs of Theorems~3 and 4 are based on the
following two lemmas. The reader should observe the subtle difference
on the right-hand sides of (\AG) and (\AGa) (the left-hand sides
being identical). 

\proclaim{Lemma 2}
For all positive integers $N$, there holds the identity
$$\multline 
\underset{\vert A\vert=N}\to{\sum _{A\subseteq[2N]} ^{}}
V(A)V\!\left(A^{-1}\right)
V(A^c)V\!\left((A^c)^{-1}\right)R\!\left(A,A^{-1}\right)
R\!\left(A^c,(A^c)^{-1}\right)
\\=
\sum _{A\subseteq[2N]} ^{}
\x^{-A}\x^{A^c}
V(A)V\!\left(A^{-1}\right)
V(A^c)V\!\left((A^c)^{-1}\right)R\!\left(A,(A^c)^{-1}\right)
R\!\left(A^c,A^{-1}\right),
\endmultline\tag\AGa$$
where $A^c$ denotes the complement of $A$ in $[2N]$, and where
$\x^{-A}$ is short for $\prod _{a\in A} ^{}x_a^{-1}$ and
$\x^{A^c}$ is short for $\prod _{a\in A^c} ^{}x_a$.
\endproclaim

\demo{Proof} 
We proceed as in the proof of Lemma~1. That is, we perform an
induction on $N$. For $N=1$, identity (\AGa) reduces to
$$\multline
2(x_1-x_1^{-1})(x_2-x_2^{-1})=
x_1x_2(x_1-x_2)(x_1^{-1}-x_2^{-1})+
x_1x_2^{-1}(x_1-x_2^{-1})(x_2-x_1^{-1})\\+
x_1^{-1}x_2(x_2-x_1^{-1})(x_1-x_2^{-1})+
x_1^{-1}x_2^{-1}(x_1^{-1}-x_2^{-1})(x_1-x_2),
\endmultline$$
which can be readily verified.
Almost all the remaining steps are identical with those in the proof of 
Lemma~1, except that more care is needed to show that, as a Laurent
polynomial in $x_1$, the degree of the right-hand side of (\AGa) is
at most $2N-1$. Indeed, by inspection, this degree is at most $2N$,
and the coefficient of $x_1^{2N}$ is equal to
$$
\sum _{A\subseteq[2N]\backslash\{1\}} ^{}
(-1)^{\vert A^c\vert-1}
V(A)V\!\left(A^{-1}\right)
V(A^c)V\!\left((A^c)^{-1}\right)R\!\left(A,(A^c)^{-1}\right)
R\!\left(A^c,A^{-1}\right),
$$
where $A^c$ now denotes the complement of $A$ in
$[2N]\backslash\{1\}$. 
(Note that only those subsets $A$ of $[2N]$ contribute to the coefficient
of $x_1^{2N}$ on the right-hand side of (\AGa) which do not contain $1$, 
whence the term
$x^{-A}x^{A^c}$ on the right-hand side of (\AGa) got cancelled due to the
contributions from the terms $R\left(A,(A^c)^{-1}\right)$ and 
$V\left((A^c)^{-1}\right)$,
respectively.)
However, in this sum, the terms indexed by $A$
respectively $A^c$ cancel each other, so that this sum does indeed
vanish.\quad \quad \qed
\enddemo

\proclaim{Lemma 3}
For all positive integers $N$, there holds the identity
$$\multline 
\underset{\vert A\vert=N}\to{\sum _{A\subseteq[2N+1]} ^{}}
V(A)V\!\left(A^{-1}\right)
V(A^c)V\!\left((A^c)^{-1}\right)R\!\left(A,A^{-1}\right)
R\!\left(A^c,(A^c)^{-1}\right)
\\=
\sum _{A\subseteq[2N+1]} ^{}
(-1)^{N+\vert A\vert}\x^{-A}\x^{A^c}
V(A)V\!\left(A^{-1}\right)
V(A^c)V\!\left((A^c)^{-1}\right)R\!\left(A,(A^c)^{-1}\right)
R\!\left(A^c,A^{-1}\right),
\endmultline\tag\AGb$$
where $A^c$ denotes the complement of $A$ in $[2N+1]$, while
$\x^{-A}$ is short for $\prod _{a\in A} ^{}x_a^{-1}$ and
$\x^{A^c}$ is short for $\prod _{a\in A^c} ^{}x_a$, as before.
\endproclaim

\demo{Proof} 
We proceed again as in the proof of Lemma~1. 
Here, the induction basis, the case $N=0$ of (\AGb), reads
$$
(x_1-x_1^{-1})=x_1-x_1^{-1}.
$$
As Laurent polynomials in $x_1$, the degree of the left-hand and
right-hand sides of (\AGb) are at most $2N+1$. This means that 
we need $4N+3$ specializations, respectively ``informations," which agree 
on both sides of (\AGb), in order to show that both sides are equal.

In the same way as this
is done in Lemma~1, one can show that, if one assumes the truth of (\AGb) 
with $N$ replaced by $N-1$, this implies that (\AGb) is true if we
specialize $x_1$ to $x_2$. Similarly to the proof of Lemma~1, it is easy to
see that both sides of (\AGb) are symmetric in the variables
$x_1,x_2,\dots,x_{2N+1}$, and also that they remain invariant up to an overall
multiplicative sign if we replace $x_1$ by $x_1^{-1}$. This means
that identity (\AGb) holds for {\it all\/} specializations of $x_1$ to
one of $x_2,x_2^{-1},x_3,x_3^{-1},\dots,x_{2N+1},x_{2N+1}^{-1}$.
These are $4N$ specializations. We still need $3$ additional
specializations, respectively ``informations," which agree on both
sides of (\AGb).

We get two more specializations by observing that both sides of (\AGb)
vanish for $x_1=\pm1$. For the left-hand side this is obvious because
of the factors $R\!\left(A,A^{-1}\right)$ respectively
$R\!\left(A^c,(A^c)^{-1}\right)$ appearing in the summand. On the
right-hand side, the summands indexed by $A$ and $A\triangle\{1\}$ 
cancel each other for $x_1=\pm1$.

Finally, the coefficients of $x_1^{2N+1}$ on the
left-hand and right-hand sides of (\AGb) are respectively
$$
\underset{\vert A\vert=N}\to{\sum _{A\subseteq[2N+1]\backslash\{1\}} ^{}}
(-1)^{N}
V(A)V\!\left(A^{-1}\right)
V(A^c)V\!\left((A^c)^{-1}\right)R\!\left(A,A^{-1}\right)
R\!\left(A^c,(A^c)^{-1}\right)
$$
and
$$
\sum _{A\subseteq[2N+1]\backslash\{1\}} ^{}
(-1)^{N}
V(A)V\!\left(A^{-1}\right)
V(A^c)V\!\left((A^c)^{-1}\right)R\!\left(A,(A^c)^{-1}\right)
R\!\left(A^c,A^{-1}\right),
$$
where, in both cases, $A^c$ now denotes the complement of $A$ in 
$[2N+1]\backslash\{1\}$. The equality of these two sums was
established in Lemma~1. This completes the proof of the 
lemma.\quad \quad \qed
\enddemo

\subhead 4. Proofs of theorems\endsubhead
This section is devoted to the proofs of Theorems~1 and 2.
The idea is to substitute the determinantal definitions (\AA)--(\ABb)
of the characters into (\AC) respectively (\ACa), expand the
determinants in the numerators by using Laplace expansion
respectively linearity of the determinant
in the rows, evaluate the resulting simpler determinants by means of
one of the Weyl denominator formulas, and reduce the
resulting expressions. By collecting appropriate terms, it is then
seen that both identities result from Lemma~1 in the preceding
section.

\demo{Proof of Theorem~1} 
We start with the left-hand side of (\AC). By (\AA), we have
$$
s_{((2m)^n)}(x_1,x_1^{-1},x_2,x_2^{-1},\dots,x_n,x_n^{-1})=
\frac {\det\limits_{1\le i,j\le 2n}
\pmatrix 
x_h^{2m\chi(t\le n)+2n-t}&1\le h\le n\\
x_{h-n}^{-(2m\chi(t\le n)+2n-t)}&n+1\le h\le 2n\\
\endpmatrix} 
{V([n])V([n]^{-1})R([n],[n]^{-1})},
\tag\AH
$$
where $\chi(\Cal S)=1$ if $\Cal S$ is
true and $\chi(\Cal S)=0$ otherwise. Here, we used the evaluation of
the Vandermonde determinant (the Weyl denominator formula for type~$A$;
cf.~\cite{\FuHaAA, p.~400 and Lemma~24.3}) in the denominator.
We now do a Laplace expansion of
the determinant along the first $n$ columns. Abbreviating the
denominator on the right-hand side of (\AH) by $D_1(n)$, this leads to
$$\multline
s_{((2m)^n)}(x_1,x_1^{-1},x_2,x_2^{-1},\dots,x_n,x_n^{-1})\\
=\frac {1} {D_1(n)}\underset{\vert A\vert+\vert B\vert=n}\to
{\sum _{A,B\subseteq[n]} ^{}}(-1)^{\big(\sum _{a\in A} a\big)+
\big(\sum _{b\in B} b\big)+n\vert
B\vert-\binom {n+1}2}\det M_1(A,B)\cdot\det M_2(A^c,B^c),
\endmultline
\tag\AI
$$
where $M_1(A,B)$ is the $n\times n$ matrix
$$\pmatrix 
x_h^{2m+2n-t}&h\in A\hfill\\
x_{h}^{-(2m+2n-t)}&h\in B\hfill\\
\endpmatrix,
$$
and $M_2(A^c,B^c)$ is the $n\times n$ matrix
$$\pmatrix 
x_h^{n-t}&h\in A^c\hfill\\
x_{h}^{-(n-t)}&h\in B^c\hfill\\
\endpmatrix,
$$
with $A^c$ denoting the complement of $A$ in $[n]$, and an analogous
meaning for $B^c$.
All determinants in (\AI) are Vandermonde determinants, except for
some trivial factors which can be taken out of the rows of
$M_1(A,B)$ respectively $M_2(A,B)$, 
and can therefore be evaluated in product form. If we
substitute the corresponding results, we obtain
$$\multline
s_{((2m)^n)}(x_1,x_1^{-1},x_2,x_2^{-1},\dots,x_n,x_n^{-1})\\
=\frac {1} {D_1(n)}\underset{\vert A\vert+\vert B\vert=n}\to
{\sum _{A,B\subseteq[n]} ^{}}(-1)^{\big(\sum _{a\in A} a\big)+
\big(\sum _{b\in B} b\big)+n\vert
B\vert-\binom {n+1}2}
\bigg(\prod _{a\in A} ^{}x_a\bigg)^{2m+n}
\bigg(\prod _{b\in B} ^{}x_b^{-1}\bigg)^{2m+n}\\
\cdot
V(A)V(B^{-1})R(A,B^{-1})V(A^c)V\!\left((B^c)^{-1}\right)
R\!\left(A^c,(B^c)^{-1}\right).
\endmultline
\tag\AJ
$$

Next we turn to the right-hand side of (\AC). By (\AB), we have
$$\align
\so_{(m^n)}(x_1,x_2,\dots,x_n)&=\frac {\det\limits_{1\le h,t\le
n}(x_h^{m+n-t+\frac 12}-x_h^{-(m+n-t+\frac 12)})}
{\det\limits_{1\le h,t\le
n}(x_h^{n-t+\frac 12}-x_h^{-(n-t+\frac 12)})}\\
&=\frac {\det\limits_{1\le h,t\le
n}(x_h^{m+n-t+\frac 12}-x_h^{-(m+n-t+\frac 12)})}
{V([n])\prod _{i=1} ^{n}x_i^{-n+i}(x_i^{1/2}-x_i^{-1/2})
\prod _{1\le i<j\le n} ^{}(x_i-x_j^{-1})},
\tag\AK
\endalign$$
where we have used the Weyl denominator formula (\AKa).
We now use linearity of the determinant in the rows. Abbreviating the
denominator on the right-hand side of (\AK) by $D_2(n)$, this leads to
$$
\so_{(m^n)}(x_1,x_2,\dots,x_n)
=\frac {1} {D_2(n)}
{\sum _{A\subseteq[n]} ^{}}(-1)^{\big(\sum _{a\in A} ^{}a\big)-\binom
{\vert A\vert+1}2}
\det M_3(A),
\tag\AL
$$
where $M_3(A)$ is the $n\times n$ matrix
$$\pmatrix 
\hphantom{-{}}x_h^{m+n-t+\frac {1} {2}}&h\in A\hfill\\
-x_{h}^{-(m+n-t+\frac {1} {2})}&h\in A^c\hfill\\
\endpmatrix,
$$
with $A^c$ denoting the complement of $A$ in $[n]$, as before.
All determinants in (\AL) are Vandermonde determinants, except for
some trivial factors which can be taken out of the rows of
$M_3(A)$, and can therefore be evaluated in product form. If we
substitute the corresponding results, we obtain
$$\multline
\so_{(m^n)}(x_1,x_2,\dots,x_n)\\
=\frac {1} {D_2(n)}
{\sum _{A\subseteq[n]} ^{}}(-1)^{\big(\sum _{a\in A} ^{}a\big)-\binom
{\vert A\vert}2+n}
\bigg(\prod _{a\in A} ^{}x_a\bigg)^{m+\frac {1} {2}}
\bigg(\prod _{a\in A^c} ^{}x_a^{-1}\bigg)^{m+\frac {1} {2}}\\
\cdot
V(A)V\left((A^c)^{-1}\right)R\!\left(A,(A^c)^{-1}\right),
\endmultline
\tag\AMa
$$
whence the product of the two characters on the right-hand side of 
(\AC) can be expanded in the form
$$\multline 
\so_{(m^n)}(x_1,x_2,\dots,x_n)
\so_{(m^n)}(-x_1,-x_2,\dots,-x_n)\\
=\frac {(-1)^{mn}} {D_1(n)}
{\sum _{A,B\subseteq[n]} ^{}}(-1)^{\big(\sum _{a\in A} ^{}a\big)+
\big(\sum _{b\in B} ^{}b\big)-\binom {\vert A\vert+1}2-\binom
{\vert B\vert}2}\kern4cm\\
\cdot
\bigg(\prod _{a\in A} ^{}x_a\bigg)^{m+\frac {1} {2}}
\bigg(\prod _{a\in A^c} ^{}x_a^{-1}\bigg)^{m+\frac {1} {2}}
\bigg(\prod _{b\in B} ^{}x_b\bigg)^{m+\frac {1} {2}}
\bigg(\prod _{b\in B^c} ^{}x_b^{-1}\bigg)^{m+\frac {1} {2}}\\
\cdot
V(A)V\left((A^c)^{-1}\right)R\!\left(A,(A^c)^{-1}\right)
V(B)V\left((B^c)^{-1}\right)R\!\left(B,(B^c)^{-1}\right).
\endmultline\tag\AM$$

We now fix disjoint subsets $A'$ and $B'$ of $[n]$ and extract the
coefficient of 
$$\bigg(\prod _{a\in A'} ^{}x_a\bigg)^{2m+n}
\bigg(\prod _{b\in B'} ^{}x_b^{-1}\bigg)^{2m+n}
\tag\AN$$
in (\AJ). (When we here say ``extract the coefficient of (\AN)," we
treat ``$m$" as if it were a formal variable. That is, in the sum on
the right-hand side of (\AJ), terms must be expanded, from each term
which results from the expansion we factor out the product (\AN), 
and, if whatever remains is independent
of $m$, it contributes to the coefficient.)
This coefficient is the subsum of the sum on the right-hand
side of (\AJ) consisting of the summands corresponding to
subsets $A$ and $B$ of $[n]$ with $A'=A\backslash B$ and 
$B'=B\backslash A$ such that $\vert A\vert+\vert B\vert=n$. We note
at this point that this implies that $C:=[n]\backslash (A'\cup B')$
must have even cardinality.
We let $A''$ be the intersection $A''=A\cap B$, so that
$A=A'\dot\cup A''$ and $B=B'\dot\cup A''$ (with $\dot\cup$ denoting
disjoint union), and we denote the complement of $A''$ in $C$
by $(A'')^c$. Since $\vert A\vert+\vert B\vert=n$, we have
$\vert A''\vert=\vert (A'')^c\vert$.
Using the above notation, then, after some manipulation, 
the above described subsum can be rewritten as
$$\multline
\frac {1} {D_1(n)}
(-1)^{n\vert B'\vert
+\frac {1} {2}\vert C\vert}
V(A')V(B')V\!\left((A')^{-1}\right)V\!\left((B')^{-1}\right)\\
\times
R\!\left(A',(B')^{-1}\right)
R\!\left(B',(A')^{-1}\right)
R\!\left(C,A'\right)
R\!\left(C,(A')^{-1}\right)
R\!\left(C^{-1},B'\right)
R\!\left(C^{-1},(B')^{-1}\right)\\
\times
\underset{\vert A''\vert=\frac {1} {2}\vert C\vert}\to
{\sum _{A''\subseteq C} ^{}}
V(A'')V\!\left((A'')^{-1}\right)
R\!\left(A'',(A'')^{-1}\right)
V\!\left((A'')^c\right)V\!\left(((A'')^c)^{-1}\right)
R\!\left((A'')^c,((A'')^c)^{-1}\right),
\endmultline
\tag\AO
$$
where we used variations of the resultant notation $R(\dots)$, 
namely
$$\align 
R\left(A,B\right)&:=\prod _{a\in A} ^{}\prod _{b\in B}
^{}(x_a-x_b),\\
R\left(A^{-1},B\right)&:=\prod _{a\in A} ^{}\prod _{b\in B}
^{}(x_a^{-1}-x_b),
\endalign$$
and
$$R\left(A^{-1},B^{-1}\right):=\prod _{a\in A} ^{}\prod _{b\in B}
^{}(x_a^{-1}-x_b^{-1}).$$

Now we turn our attention to (\AM). First of all, we observe that
by exchanging the roles of $A$ and $B$ in the summand on the
right-hand side of (\AM), the summand does not change except for a
sign of $(-1)^{\vert A\vert-\vert B\vert}$. Thus, all summands
corresponding to sets $A$ and $B$ whose cardinalities do not have the
same parity cancel each other. We can therefore restrict the sum in
(\AM) to subsets $A$ and $B$ of $[n]$ with $\vert A\vert\equiv\vert
B\vert$~(mod~$2$).
Consequently, if we extract the coefficient of (\AN) in (\AM), 
then, after some manipulation, we obtain
$$\multline
\frac {1} {D_1(n)}
(-1)^{mn+n\vert B'\vert
+\frac {1} {2}\vert C\vert}
V(A')V(B')V\!\left((A')^{-1}\right)V\!\left((B')^{-1}\right)\\
\times
R\!\left(A',(B')^{-1}\right)
R\!\left(B',(A')^{-1}\right)
R\!\left(C,A'\right)
R\!\left(C,(A')^{-1}\right)
R\!\left(C^{-1},B'\right)
R\!\left(C^{-1},(B')^{-1}\right)\\
\times
{\sum _{A''\subseteq C} ^{}}
V(A'')V\!\left((A'')^{-1}\right)
R\!\left(A'',((A'')^c)^{-1}\right)
V\!\left((A'')^c\right)V\!\left(((A'')^c)^{-1}\right)
R\!\left((A'')^c,(A'')^{-1}\right).
\endmultline
\tag\AP
$$
Here, with $A$ and $B$ as in (\AM), the meanings of $A'$, $B'$,
$A''$,
and $C$ are $A'=A\cap B$, $B'=[n]\backslash (A\cup B)$, 
$A''=A\backslash B$, 
$C=[n]\backslash A'\cup B'$, and $(A'')^c$ is the complement of $A''$ in
$C$, so that $A=A'\dot\cup A''$ and
$B=A'\dot\cup (A'')^c$. Because of the restriction
$\vert A\vert\equiv\vert B\vert$~(mod~$2$), again $C$ must
have even cardinality.
Clearly, clearing the factors common to (\AO) and (\AP), 
the equality of (\AO) and (\AP) is a direct consequence of
(\AG). This completes the proof of the theorem.\quad \quad
\qed
\enddemo

\demo{Proof of Theorem~2} 
As we have shown when we rewrote Theorems~1 and 2 uniformly as
(\ACd), Theorem~2 results from the proof of Theorem~1 by ``replacing
$m$ by $m+\frac {1} {2}$." Since, in the proof of Theorem~1, it is
nowhere used that $m$ is an integer, and since, in fact, $m$ is
treated there as a formal variable, Theorem~2 follows immediately.\quad 
\quad \qed
\enddemo

\subhead 5. Combinatorial interpretations\endsubhead
The purpose of this section is to give combinatorial
interpretations of Theorems~1 and 2. As we already said in the
introduction, these interpretations were at the origin of this work,
which suggested Theorems~1 and 2 in the first place.
They involve {\it rhombus tilings}, respectively {\it plane partitions}. 
When we speak of a
rhombus tiling of some region, we always mean 
a tiling of the region by unit rhombi with angles of
$60^\circ$ and $120^\circ$. Examples of such tilings can be found in
Figures~\FA--\FE. (Dotted lines and shadings should be ignored at this point.) 
We shall not recall the relevant plane partition definitions here,
but instead refer the reader to \cite{\BresAO}, and to \cite{\KupeAA}
for explanations on the relation between plane partitions and rhombus
tilings of hexagons. 

We claim that,
by specializing $x_1=x_2=\dots=x_n=1$ in (\AC), we obtain the
combinatorial factorization
$$
\PP(2m,n,n)=
\SPP(2m,n,n)\cdot
\TCPP(2m,n,n),
\tag\AD$$
where $\PP(2m,n,n)$ denotes the number of plane partitions contained 
in the $(2m)\times n\times n$ box (or, equivalently, the number of
{\it rhombus tilings} of a hexagon with side lengths $2m,n,n,2m,n,n$;
see Figure~\FA\ for an example in which $m=2$ and $n=3$),
$\SPP(2m,n,n)$ denotes the number of {\it symmetric plane partitions} contained 
in the $(2m)\times n\times n$ box (or, equivalently, the number of
rhombus tilings of a hexagon with side lengths $2m,n,n,2m,n,n$ which
are symmetric with respect to the vertical symmetry axis of the
hexagon; see Figure~\FB\ for an example in which $m=2$ and $n=3$),
and $\TCPP(2m,n,n)$ denotes the number of {\it transpose complementary
plane partitions} contained 
in the $(2m)\times n\times n$ box (or, equivalently, the number of
rhombus tilings of a hexagon with side lengths $2m,n,n,2m,n,n$ which
are symmetric with respect to the horizontal symmetry axis of the
hexagon; see Figure~\FC.a for an example in which $m=n=3$; the dotted lines
should be ignored for the moment).\footnote{Clearly, the
factorization (\AD) could be readily verified directly by using the
known product formulas for $\PP(2m,n,n)$, $\SPP(2m,n,n)$, and
$\TCPP(2m,n,n)$ (cf\. \cite{\BresAO}). However, the point here is
that it is a consequence of the more general factorization (\AC)
featuring a Schur function and odd orthogonal characters.}

\midinsert
\vskip10pt
\vbox{\noindent
\centertexdraw{
\drawdim truecm \setunitscale.7
\linewd.15
\move(0 4)
\RhombusB \RhombusA \RhombusB \RhombusB \RhombusA \RhombusB \RhombusA 
\move(0.866025 4.5)
\RhombusA \RhombusB \RhombusB \RhombusA \RhombusB \RhombusB \RhombusA 
\move(1.73205 5)
\RhombusA \RhombusB \RhombusB \RhombusA \RhombusB \RhombusA \RhombusB
\move(0 3)
\RhombusC 
\move(0 2)
\RhombusC 
\move(0 1)
\RhombusC \RhombusC 
\move(3.4641 5)
\RhombusC \RhombusC 
\move(3.4641 4)
\RhombusC \RhombusC 
\move(4.33013 2.5)
\RhombusC 
\rtext td:60 (4.1 4.6){$\sideset {} \and {} \to 
    {\left.\vbox{\vskip1.1cm}\right\}}$}
\rtext td:120 (1.34 4.75){$\sideset {} \and {} \to 
    {\left.\vbox{\vskip1.1cm}\right\}}$}
\rtext td:0 (-1.5 1.8){$\sideset {2m} \and {}\to 
    {\left\{\vbox{\vskip1.3cm}\right.}$}
\htext (.55 5.2){$n$}
\htext (4.25 5.2){$n$}
}
\centerline{\eightpoint A rhombus tiling of a hexagon}
\vskip8pt
\centerline{\eightpoint Figure \FA}
}
\vskip10pt
\endinsert

We start with the left-hand side of (\AD).
The fact that the specialized Schur
function $s_{((2m)^n)}(1,1,\dots,1)$ (with $2n$ occurrences of $1$ in
the argument) is equal to the number of plane partitions in the 
$(2m)\times n\times n$ box, and that this is equal to the number of
rhombus tilings of a hexagon with side lengths $2m,n,n,2m,n,n$, is
well-known (cf\. \cite{\KupeAA} and \cite{\StanBI, Sec.~7.21}).

\midinsert
\vskip10pt
\vbox{\noindent
\centertexdraw{
\drawdim truecm \setunitscale.7
\linewd.15
\move(0 4)
\RhombusB \RhombusB \RhombusB \RhombusA \RhombusB \RhombusA \RhombusA
\move(0.866025 4.5)
\RhombusB \RhombusA \RhombusB \RhombusB \RhombusA \RhombusB \RhombusA 
\move(1.73205 5)
\RhombusB \RhombusA \RhombusA \RhombusB \RhombusB \RhombusA \RhombusB
\move(0 1)
\RhombusC 
\move(0.866025 3.5)
\RhombusC 
\move(2.59808 3.5)
\RhombusC 
\move(4.33013 3.5)
\RhombusC 
\move(4.33013 2.5)
\RhombusC 
\move(2.59808 5.5)
\RhombusC \RhombusC \RhombusC 
\rtext td:60 (4.1 4.6){$\sideset {} \and {} \to 
    {\left.\vbox{\vskip1.1cm}\right\}}$}
\rtext td:120 (1.34 4.75){$\sideset {} \and {} \to 
    {\left.\vbox{\vskip1.1cm}\right\}}$}
\rtext td:0 (-1.5 1.8){$\sideset {2m} \and {}\to 
    {\left\{\vbox{\vskip1.3cm}\right.}$}
\htext (.55 5.2){$n$}
\htext (4.25 5.2){$n$}
}
\centerline{\eightpoint A vertically symmetric rhombus tiling of a hexagon}
\vskip8pt
\centerline{\eightpoint Figure \FB}
}
\vskip10pt
\endinsert

It is also well-known (see \cite{\BresAO, Sec.~4.3} or 
\cite{\MacdAC, Ch.~I, Sec.~5, Ex.~15--17}) that the number of symmetric plane partitions contained 
in the $(2m)\times n\times n$ box is equal to 
$\so_{(m^n)}(1,1,\dots,1)$ (with $n$ occurrences of $1$ in
the argument). 

The argument which explains that $\so_{(m^n)}(-1,-1,\dots,-1)$
(with $n$ occurrences of $-1$ in the argument) counts 
transpose complementary plane partitions contained 
in the $(2m)\times n\times n$ box (up to sign) is more elaborate. 
We begin by
setting $x_i=-q^{i-1}$ in (\AB) with $N=n$ and $\la=(m^n)$, to obtain
$$
\so_{(m^n)}(-q,-q^2,\dots,-q^n)=(-1)^{mn}
\frac {\det\limits_{1\le h,t\le
n}(q^{(h-1)(m+n-t+\frac 12)}+q^{-(h-1)(m+n-t+\frac 12)})}
{\det\limits_{1\le h,t\le
n}(q^{(h-1)(n-t+\frac 12)}+q^{-(h-1)(n-t+\frac 12)})}.
$$
The determinants on the right-hand side can be evaluated by means of
(\AR). (This is seen by replacing $h$ by $n+1-h$ and then
interchanging the roles of $h$ and $t$ in the determinants above.)
If we subsequently let $q$ tend to $1$, then we obtain
$$
\so_{(m^n)}(-1,-1,\dots,-1)=(-1)^{mn}
\prod _{1\le h<t\le n} ^{}\frac {2m+2n+1-h-t} {2n+1-h-t}.
$$
On the other hand, the product on the right-hand side is equal to a
specialized symplectic character. Namely, 
by specializing $\la=(m^N)$, $x_i=q^i$, $i=1,2,\dots,N$, in (\ABa),
using the identity (\AS),
and finally letting $q$ tend to $1$, it is seen that
$$
\sp_{(m^N)}(1,1,\dots,1)=
\prod _{1\le h<t\le N+1} ^{}\frac {2m+2N+3-h-t} {2N+3-h-t}
$$
(with $N$ occurrences of $1$ in the argument). Thus, for $N=n-1$, we
obtain that
$$\so_{(m^n)}(-1,-1,\dots,-1)=(-1)^{mn}
\sp_{(m^{n-1})}(1,1,\dots,1).$$

\midinsert
\vskip10pt
\vbox{\noindent
\centertexdraw{
\drawdim truecm \setunitscale.7
\linewd.15
\move(0 4)
\RhombusB \RhombusB \RhombusB \RhombusA \RhombusB \RhombusA \RhombusB
\RhombusA \RhombusB
\move(0.866025 4.5)
\RhombusB \RhombusB \RhombusA \RhombusB \RhombusA \RhombusB \RhombusB
\RhombusA \RhombusB
\move(1.73205 5)
\RhombusB \RhombusA \RhombusA \RhombusB \RhombusB \RhombusA \RhombusB
\RhombusB \RhombusB
\move(0 1)
\RhombusC 
\move(0 0)
\RhombusC \RhombusC 
\move(0 -1)
\RhombusC \RhombusC \RhombusC 
\move(3.4641 1)
\RhombusC 
\move(1.73205 4)
\RhombusC \RhombusC 
\move(2.59808 5.5)
\RhombusC \RhombusC \RhombusC 
\move(4.33013 3.5)
\RhombusC 
\move(4.33013 2.5)
\RhombusC 
\ringerl(0 3.5)
\ringerl(0 2.5)
\ringerl(0 1.5)
%
\ringerl(5.19615 3.5)
\ringerl(5.19615 2.5)
\ringerl(5.19615 1.5)
%
\linewd.08
\move(0 3.5)
\odaSchritt \odaSchritt \odaSchritt \hdSchritt \hdSchritt \hdSchritt
\move(0 2.5)
\odaSchritt \odaSchritt \hdSchritt \hdSchritt \odaSchritt \hdSchritt
\move(0 1.5)
\odaSchritt \hdSchritt \odaSchritt \hdSchritt \odaSchritt \hdSchritt
%
%
\rtext td:60 (4.1 4.6){$\sideset {} \and {} \to 
    {\left.\vbox{\vskip1.1cm}\right\}}$}
\rtext td:120 (1.34 4.75){$\sideset {} \and {} \to 
    {\left.\vbox{\vskip1.1cm}\right\}}$}
\rtext td:0 (-1.45 .9){$\sideset {2m} \and {}\to 
    {\left\{\vbox{\vskip1.85cm}\right.}$}
\htext (.55 5.2){$n$}
\htext (4.25 5.2){$n$}
\move(9 0)
\bsegment
\ringerl(0 3.5)
\ringerl(0 2.5)
\ringerl(0 1.5)
%
\ringerl(5.19615 3.5)
\ringerl(5.19615 2.5)
\ringerl(5.19615 1.5)
%
\move(0 3.5)
\odaSchritt \odaSchritt \odaSchritt \hdSchritt \hdSchritt \hdSchritt
\move(0 2.5)
\odaSchritt \odaSchritt \hdSchritt \hdSchritt \odaSchritt \hdSchritt
\move(0 1.5)
\odaSchritt \hdSchritt \odaSchritt \hdSchritt \odaSchritt \hdSchritt
\esegment
}
\centerline{\eightpoint a. A horizontally symmetric rhombus tiling of a 
hexagon\quad 
b. The corresponding path family}
$$
\Gitter(4,7)(-4,0)
\Koordinatenachsen(4,7)(-4,0)
\Pfad(-1,0),212121\endPfad
\Pfad(-2,1),221121\endPfad
\Pfad(-3,2),222111\endPfad
\DickPunkt(-1,0)
\DickPunkt(-2,1)
\DickPunkt(-3,2)
\DickPunkt(2,3)
\DickPunkt(1,4)
\DickPunkt(0,5)
\Label\ro{P_1}(-1,1)
\Label\ro{P_2}(-1,3)
\Label\ro{P_3}(-2,5)
\thinlines
\Diagonale(-2,-1)5
\hbox{\hskip6cm}
\Gitter(4,7)(-4,0)
\Koordinatenachsen(4,7)(-4,0)
\Pfad(-1,1),1212\endPfad
\Pfad(-2,2),2112\endPfad
\Pfad(-3,3),2211\endPfad
\DickPunkt(-1,1)
\DickPunkt(-2,2)
\DickPunkt(-3,3)
\DickPunkt(1,3)
\DickPunkt(0,4)
\DickPunkt(-1,5)
\Label\ro{P'_1}(-1,1)
\Label\ro{P'_2}(-1,3)
\Label\ro{P'_3}(-2,5)
\thinlines
\Diagonale(-2,-1)5
$$
\centerline{\eightpoint c. Orthogonal lattice paths
\hskip2cm
d. Shorter non-intersecting lattice paths}
\vskip8pt
\centerline{\eightpoint Figure \FC}
}
\vskip10pt
\endinsert

To conclude the argument, it is easy to see using families of $m$
lattice paths starting at the mid-points of the
top-most $m$ edges along the left vertical side of length $2m$ 
of the hexagon with side lengths $2m,n,n,2m,n,n$ and ending 
at the mid-points of the top-most $m$ edges along the right 
vertical side of length $2m$ of the hexagon, the paths ``following
the rhombi" of the tiling (and, thus, staying necessarily in
the upper half of the hexagon; see Figure~\FC.a,b), that rhombus tilings 
of a hexagon with side lengths $2m,n,n,2m,n,n$ which
are symmetric with respect to the horizontal symmetry axis of the
hexagon are in bijection with families $(P_1,P_2,\dots,P_m)$ of
non-intersecting lattice paths consisting of horizontal and vertical
unit steps, the path $P_i$ starting at $(-i,i-1)$
and ending at $(n-i,n+i-1)$, $i=1,2,\dots,m$, all paths never passing
below the line $y=x+1$ (see Figure~\FC.c for the path family
resulting from the one in Figure~\FC.b by deforming the paths to
orthogonal ones; as usual, the term
``non-intersecting" means that no two paths in the family have a
common point). Because of the boundary $y=x+1$ and the
condition that paths are non-intersecting, the initial step of {\it
any} path in such a family must be a vertical step while the final
step of {\it any} path must be a horizontal step. Thus, our rhombus
tilings are in bijection with families $(P_1,P_2,\dots,P_m)$ of
non-intersecting lattice paths consisting of horizontal and vertical
unit steps, the path $P_i$ starting at $(-i,i)$
and ending at $(n-1-i,n-1+i)$, $i=1,2,\dots,m$, all paths never passing
below the line $y=x+1$ (see Figure~\FC.d).
It is known (see \cite{\FuKrAA, Sec.~5}) that these families of
non-intersecting lattice paths are counted by 
$\sp_{(m^{n-1})}(1,1,\dots,1)$ (with $n$ occurrences of $1$ in the
argument).
Thus, $\vert\so_{(m^n)}(-1,-1,\dots,-1)\vert$ 
counts indeed rhombus tilings 
of a hexagon with side lengths $2m,n,n,2m,n,n$ which
are symmetric with respect to the horizontal symmetry axis of the
hexagon, respectively transpose
complementary plane partitions in the $(2m)\times n\times n$ box.

\medskip
For Theorem~2 we are also able to provide a combinatorial
interpretation in the context of rhombus tilings. 
However, it may be less convincing. 

We specialize $x_1=x_2=\dots=x_n=1$ in (\ACa). Clearly, 
$s_{((2m+1)^n)}(1,1,\dots,1)$ (with $2n$ occurrences of $1$ in
the argument) counts plane partitions in the 
$(2m+1)\times n\times n$ box, respectively
rhombus tilings of a hexagon with side lengths $2m+1,n,n,2m+1,n,n$.

\midinsert
\vskip10pt
\vbox{\noindent
\centertexdraw{
\drawdim truecm \setunitscale.7
\linewd.15
\move(0 4)
\RhombusA \RhombusB \RhombusB \RhombusA \RhombusA \RhombusB 
\RhombusB \RhombusA
\move(0.866025 4.5)
\RhombusA \RhombusB \RhombusB \RhombusA \RhombusA \RhombusA
\RhombusB \RhombusB
\move(1.73205 5)
\RhombusA \RhombusA \RhombusB \RhombusB \RhombusA 
\move(0 1)
\RhombusB \RhombusB \RhombusA \RhombusA \RhombusA
\move(0 4)
\RhombusC 
\move(0 3)
\RhombusC 
\move(0.866025 1.5)
\RhombusC \RhombusC 
\move(0.866025 .5)
\RhombusC \RhombusC 
\move(3.4641 1)
\RhombusC 
\move(2.59808 4.5)
\RhombusC 
\move(4.33013 4.5)
\RhombusC 
\move(4.33013 3.5)
\RhombusC 
\rtext td:60 (4.1 4.6){$\sideset {} \and {} \to 
    {\left.\vbox{\vskip1.1cm}\right\}}$}
\rtext td:120 (1.34 4.75){$\sideset {} \and {} \to 
    {\left.\vbox{\vskip1.1cm}\right\}}$}
\rtext td:0 (-2.45 1.3){$\sideset {2m+1} \and {}\to 
    {\left\{\vbox{\vskip1.6cm}\right.}$}
\htext (.55 5.2){$n$}
\htext (4.25 5.2){$n$}
}
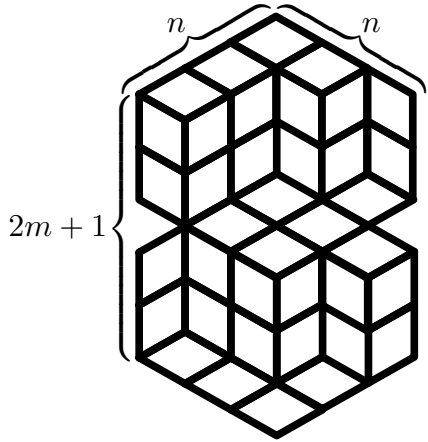
\centerline{\eightpoint A horizontally symmetric rhombus tiling of a hexagon
with two missing triangles}
\vskip8pt
\centerline{\eightpoint Figure \FD}
}
\vskip10pt
\endinsert

Moreover, the argument above shows that $\sp_{(m^n)}(1,1,\dots,1)$
(with $n$ occurrences of $1$ in the argument) counts rhombus tilings
of a hexagon with side lengths $2m+1,n,n,2m+1,n,n$ from which two unit
triangles have been removed on the left and the right end of the
horizontal symmetry axis of the hexagon, the tilings being symmetric
with respect to this axis 
(see Figure~\FD\ for an example in which $m=2$ and $n=3$).
Alternatively, this is also the number of horizontally symmetric rhombus 
tilings of a (full) hexagon with side lengths $2m,n+1,n+1,2m,n+1,n+1$ 
(compare with Figure~\FC.a).

In order to find a combinatorial interpretation of 
$o^{even}_{((m+1)^n)}(1,1,\dots,1)$ (with $n$ occurrences of $1$ in the 
argument), we start with the decomposition
$$o^{even}_{((m+1)^n)}(1,1,\dots,1)=
\so^{even}_{((m+1)^n)}(1,1,\dots,1)
+\so^{even}_{((m+1)^{n-1},-m-1)}(1,1,\dots,1),$$
where
$$
\so^{even}_\la(x_1,x_2,\dots,x_n)
=\frac {\det\limits_{1\le h,t\le
n}(x_h^{\la_t+n-t}+x_h^{-(\la_t+n-t)})+
\det\limits_{1\le h,t\le
n}(x_h^{\la_t+n-t}-x_h^{-(\la_t+n-t)})} 
{\det\limits_{1\le h,t\le
n}(x_h^{n-t}+x_h^{-(n-t)})}
$$
is an irreducible character of $SO_{2n}(\C)$ (and its spin covering
group; see \cite{\FuHaAA, (24.40)}). 
Here, $\la=(\la_1,\la_2,\dots,\la_n)$ is a non-increasing sequence 
of (possibly negative) integers or half-integers with 
$\la_{n-1}\ge\vert\la_n\vert$.
It was proved in \cite{\BrGrAA} (see \cite{\KratBC} for a
common generalization) that, for all non-negative integers or
half-integers $c$, we have
$$\so^{even}_{(c^n)}(x_1,x_2,\dots,x_n)=
(x_1x_2\cdots x_n)^{-c}\cdot\underset
\oddcols\!\big(((2c)^n)/\nu\big)=0
\to{\sum _{\nu\subseteq ((2c)^n)} ^{}}s_\nu(x_1,x_2,\dots,x_n)
\tag\AE$$
and
$$\so^{even}_{(c^{n-1},-c)}(x_1,x_2,\dots,x_n)=
(x_1x_2\cdots x_n)^{-c}\cdot\underset
\oddcols\!\big(((2c)^n)/\nu\big)=2c
\to{\sum _{\nu\subseteq ((2c)^n)} ^{}}s_\nu(x_1,x_2,\dots,x_n),
\tag\AF$$
where $\oddcols\!\big(((2c)^n)/\nu\big)$ denotes the number of odd
columns of the {\it skew diagram} (cf\. \cite{\MacdAC, p.~4})
$((2c)^n)/\nu$. Thus, we have
$$o^{even}_{((m+1)^n)}(1,1,\dots,1)=
\sum _{\nu\subseteq ((2m+2)^n)} ^{}\kern-15pt{}^{\displaystyle\prime}
\kern15pt
s_\nu(1,1,\dots,1),
$$
where $\sum{}^{\textstyle\prime}$ is taken over all diagrams $\nu$
with the property that $((2m+2)^n)/\nu$ consists either of only even
columns or of only odd columns. (Both the even orthogonal character
on the left-hand side and the Schur functions on the right-hand
side contain $n$ occurrences of $1$ in their arguments.)

\midinsert
\vbox{\noindent
\hbox{\hskip-3.5cm
\centertexdraw{
\drawdim truecm \setunitscale.7
\linewd.15
\move(0 4)
\RhombusB \RhombusA \RhombusB \RhombusA \RhombusB 
\RhombusB \RhombusA \RhombusA \RhombusB \RhombusB 
\move(0.866025 4.5)
\RhombusB \RhombusA \RhombusA \RhombusB \RhombusB 
\RhombusB \RhombusA \RhombusA \RhombusB \RhombusB 
\move(1.73205 5)
\RhombusB \RhombusA \RhombusA \RhombusA \RhombusB 
\RhombusA \RhombusB \RhombusB \RhombusB \RhombusB 
\move(2.59808 5.5)
\RhombusB \RhombusA \RhombusA \RhombusA \RhombusA
\RhombusB \RhombusB \RhombusB \RhombusB \RhombusB 
\move(0 3)
\RhombusC 
\move(0 2)
\RhombusC \RhombusC 
\move(0 1)
\RhombusC \RhombusC 
\move(0 0)
\RhombusC \RhombusC \RhombusC \RhombusC 
\move(0 -1)
\RhombusC \RhombusC \RhombusC \RhombusC 
\move(3.4641 3)
\RhombusC 
\move(3.4641 2)
\RhombusC \RhombusC 
\move(3.4641 1)
\RhombusC \RhombusC 
\move(3.4641 6)
\RhombusC \RhombusC \RhombusC \RhombusC 
\ringerl(0 3.5)
\ringerl(0 2.5)
\ringerl(0 1.5)
\ringerl(0 .5)
\ringerl(0 -.5)
\ringerl(0 -1.5)
\ringerl(3.4641 5.5)
\ringerl(3.4641 2.5)
\ringerl(3.4641 1.5)
\ringerl(3.4641 .5)
\ringerl(3.4641 -2.5)
\ringerl(3.4641 -3.5)
\linewd.08
\move(0 3.5)
\odaSchritt \odaSchritt \odaSchritt \odaSchritt 
\move(0 2.5)
\hdSchritt \odaSchritt \hdSchritt \odaSchritt 
\move(0 1.5)
\hdSchritt \hdSchritt \odaSchritt \odaSchritt
\move(0 .5)
\hdSchritt \hdSchritt \odaSchritt \odaSchritt
\move(0 -.5)
\hdSchritt \hdSchritt \hdSchritt \hdSchritt 
\move(0 -1.5)
\hdSchritt \hdSchritt \hdSchritt \hdSchritt 
%
\move(3.4641 5)  \rlvec(-.866025 -.5) \rlvec(.866025 -.5)
\lfill f:.3
\rlvec(.866025 .5) \rlvec(-.866025 .5)
\lfill f:.3
\move(3.4641 4)  \rlvec(-.866025 -.5) \rlvec(.866025 -.5)
\lfill f:.3
\rlvec(.866025 .5) \rlvec(-.866025 .5)
\lfill f:.3
\move(3.4641 0)  \rlvec(-.866025 -.5) \rlvec(.866025 -.5)
\lfill f:.3
\rlvec(.866025 .5) \rlvec(-.866025 .5)
\lfill f:.3
\move(3.4641 -1)  \rlvec(-.866025 -.5) \rlvec(.866025 -.5)
\lfill f:.3
\rlvec(.866025 .5) \rlvec(-.866025 .5)
\lfill f:.3
\rtext td:60 (5.3 4.9){$\sideset {} \and {} \to 
    {\left.\vbox{\vskip1.3cm}\right\}}$}
\rtext td:120 (1.94 5.1){$\sideset {} \and {} \to 
    {\left.\vbox{\vskip1.3cm}\right\}}$}
\rtext td:0 (-2.5 .9){$\sideset {2m+2} \and {}\to 
    {\left\{\vbox{\vskip1.9cm}\right.}$}
\htext (1.15 5.55){$n$}
\htext (5.25 5.55){$n$}
}
\hskip-4cm\hbox {}}

\vskip-6cm
$$
\hbox{\hskip7.5cm}
\Gitter(4,11)(-6,0)
\Koordinatenachsen(4,11)(-6,0)
\Pfad(-1,1),1111\endPfad
\Pfad(-2,2),1111\endPfad
\Pfad(-3,3),1122\endPfad
\Pfad(-4,4),1122\endPfad
\Pfad(-5,5),1212\endPfad
\Pfad(-6,6),2222\endPfad
\DickPunkt(-1,1)
\DickPunkt(-2,2)
\DickPunkt(-3,3)
\DickPunkt(-4,4)
\DickPunkt(-5,5)
\DickPunkt(-6,6)
\DickPunkt(3,1)
\DickPunkt(2,2)
\DickPunkt(-1,5)
\DickPunkt(-2,6)
\DickPunkt(-3,7)
\DickPunkt(-6,10)
\Label\o{P_1}(1,0)
\Label\ro{P_2}(0,2)
\Label\ro{P_3}(-1,4)
\Label\ro{P_4}(-3,4)
\Label\ro{P_5}(-4,6)
\Label\ro{P_6}(-6,7)
$$
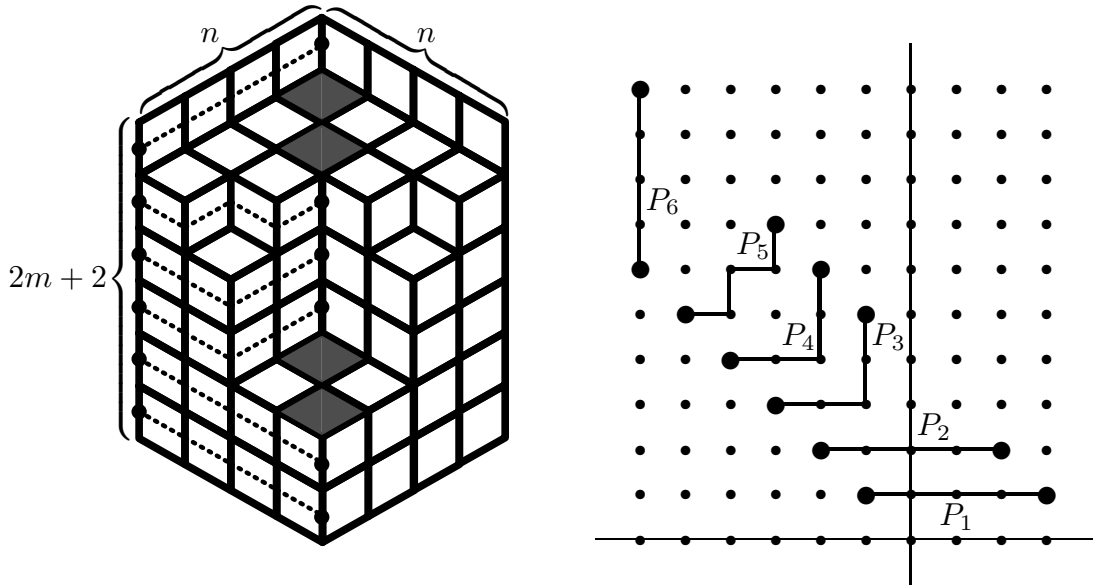
\centerline{\eightpoint a. A vertically symmetric rhombus tiling of a
hexagon
\hskip1cm
b. Non-intersecting lattice paths\hskip1cm}
\vskip8pt
\centerline{\eightpoint Figure \FE}
}
\vskip10pt
\endinsert

Given a shape $\nu$ contained in $((2m+2)^n)$,
it is well-known (cf\. \cite{\FuKrAA, Sec.~4} or \cite{\SagaAL, Ch.~4})
that $s_\nu(1,1,\dots,1)$ counts families $(P_1,\dots,P_{2m+2})$ of 
non-intersecting lattice paths consisting of horizontal and vertical
unit steps, where the path $P_i$ runs from $(-i,i)$ to
$(\nu_i'-i,n-\nu_i'+i)$, $i=1,2,\dots,2m+2$. (Here, $\nu'$
denotes the partition {\it conjugate} to $\nu$; cf\. \cite{\MacdAC,
p.~2}). On the other hand, in a similar way as above,
such families of non-intersecting lattice 
paths (if $\nu$ is allowed to be any partition) 
are in bijection with rhombus tilings of a hexagon 
with side lengths $2m+2,n,n,2m+2,n,n$ which are symmetric with
respect to the vertical symmetry axis of the hexagon (see
Figure~\FE\ for an example in which $m=2$, $n=4$, $\nu=(5,5,2,2)$).
The property that $((2m+2)^n)/\nu$ consists either of only even
columns or of only odd columns translates into the property that, in the
corresponding rhombus tilings, chains of successive
horizontally oriented rhombi along the vertical symmetry axis of the
hexagon can appear in the interior of the hexagon only if they have
even length. (By definition, a ``chain of horizontally oriented
rhombi along the vertical symmetry axis of the hexagon" is a set of
horizontally oriented rhombi sitting on the vertical symmetry axis
which, together, form a topologically connected set. For a chain, to be in the
interior, means that none of the rhombi of the chain touches the boundary of 
the hexagon.
In Figure~\FE.a, there are two such chains. They consist of
the two contiguous strings of grey shaded rhombi, respectively.) 

In summary, Theorem~2, when specialized to $x_1=x_2=\dots=x_n=1$, can
be interpreted combinatorially as follows: the term
$s_{((2m+1)^n)}(1,1,\dots,1)$ on the left-hand side counts
rhombus tilings of a hexagon with side lengths $2m+1,n,n,2m+1,n,n$.
The term $\sp_{(m^n)}(1,1,\dots,1)$ counts horizontally symmetric
rhombus tilings
of a hexagon with side lengths $2m+1,n,n,2m+1,n,n$ from which two unit
triangles have been removed on the left and the right end of the
horizontal symmetry axis of the hexagon. Finally, the term
$o^{even}_{((m+1)^n)}(1,1,\dots,1)$ counts vertically symmetric
rhombus tilings of a hexagon with side lengths $2m+2,n,n,2m+2,n,n$ 
with the property that chains of successive
horizontally oriented rhombi along the vertical symmetry axis of the
hexagon can appear in the interior of the hexagon only if they have
even length. We remark that these rhombus tilings could equivalently be seen as
symmetric plane partitions in a $(2m+2)\times n\times n$ box in
which ``central terraces" (that is, horizontal levels situated along
the plane of symmetry) must have even length, except at height $0$
and at height $2m+2$. In other words, this specialization yields the identity
$$
\PP(2m+1,n,n) = \TCPP(2m,n+1,n+1) \SPP^*(2m+2,n,n),
\tag\AFb
$$
where the $^*$ indicates that only those symmetric plane partitions are counted which
satisfy the ``even central terraces in the interior" condition described above.
In particular, by specializing $x_h=q^{h-1}$ in
(\ABb), using (\AR) for evaluating the two determinants in (\ABb), and
then letting $q$ tend to $1$, we obtain that
$$\SPP^*(2m,n,n)=2\prod _{1\le h<t\le n} ^{}\frac {2m+2n-h-t} {2n-h-t}.
\tag\AFa$$

\subhead 6. More factorization theorems\endsubhead
In this final section, we present two further factorization theorems
similar to those of Theorems~1 and 2, which have been hinted at to us by
Ron King. Since these theorems can be proved in a manner very 
similar to the one which led to our proofs of Theorems~1 and 2, we
content ourselves with giving only sketches of proofs, leaving the
(easily filled in) details to the reader.

\proclaim{Theorem 3}For any non-negative integers $m$ and $n$, we have
$$\multline 
s_{((2m+1)^n)}(x_1,x_1^{-1},x_2,x_2^{-1},\dots,x_n,x_n^{-1})+
s_{((2m+1)^{n-1})}(x_1,x_1^{-1},x_2,x_2^{-1},\dots,x_n,x_n^{-1})\\
=(-1)^{mn}
\so_{((m+1)^n)}(x_1,x_2,\dots,x_n)\,
\so_{(m^n)}(-x_1,-x_2,\dots,-x_n).
\endmultline
\tag\AT$$
\endproclaim

\demo{Sketch of proof} 
We proceed in the same way as in the proof of Theorem~1. By
comparison with (\AJ), we see that the left-hand side of (\AT)
is equal to
$$\multline
\frac {1} {D_1(n)}\underset{\vert A\vert+\vert B\vert=n}\to
{\sum _{A,B\subseteq[n]} ^{}}(-1)^{\big(\sum _{a\in A} a\big)+
\big(\sum _{b\in B} b\big)+n\vert
B\vert-\binom {n+1}2}
\bigg(\prod _{a\in A} ^{}x_a\bigg)^{2m+n+1}
\bigg(\prod _{b\in B} ^{}x_b^{-1}\bigg)^{2m+n+1}\\
\cdot
V(A)V(B^{-1})R(A,B^{-1})V(A^c)V\!\left((B^c)^{-1}\right)
R\!\left(A^c,(B^c)^{-1}\right)\\\kern-1pt
+
\frac {1} {D_1(n)}\underset{\vert A\vert+\vert B\vert=n-1}\to
{\sum _{A,B\subseteq[n]} ^{}}(-1)^{\big(\sum _{a\in A} a\big)+
\big(\sum _{b\in B} b\big)+n\vert
B\vert-\binom {n+1}2}
\bigg(\prod _{a\in A} ^{}x_a\bigg)^{2m+n+2}
\bigg(\prod _{b\in B} ^{}x_b^{-1}\bigg)^{2m+n+2}\\
\cdot
V(A)V(B^{-1})R(A,B^{-1})V(A^c)V\!\left((B^c)^{-1}\right)
R\!\left(A^c,(B^c)^{-1}\right),
\endmultline\tag\AU
$$
with $D_1(n)$ having the same meaning as in the proof of Theorem~1.
On the other hand, from (\AMa) we see that the right-hand side of
(\AT) is equal to
$$\multline 
\frac {(-1)^{mn}} {D_1(n)}
{\sum _{A,B\subseteq[n]} ^{}}(-1)^{\big(\sum _{a\in A} ^{}a\big)+
\big(\sum _{b\in B} ^{}b\big)-\binom {\vert A\vert+1}2-\binom
{\vert B\vert}2}\kern4cm\\
\cdot
\bigg(\prod _{a\in A} ^{}x_a\bigg)^{m+\frac {1} {2}}
\bigg(\prod _{a\in A^c} ^{}x_a^{-1}\bigg)^{m+\frac {1} {2}}
\bigg(\prod _{b\in B} ^{}x_b\bigg)^{m+\frac {3} {2}}
\bigg(\prod _{b\in B^c} ^{}x_b^{-1}\bigg)^{m+\frac {3} {2}}\\
\cdot
V(A)V\left((A^c)^{-1}\right)R\!\left(A,(A^c)^{-1}\right)
V(B)V\left((B^c)^{-1}\right)R\!\left(B,(B^c)^{-1}\right)\\
=\frac {(-1)^{mn}} {D_1(n)}
{\sum _{A,B\subseteq[n]} ^{}}(-1)^{\big(\sum _{a\in A} ^{}a\big)+
\big(\sum _{b\in B} ^{}b\big)-\binom {\vert A\vert+1}2-\binom
{\vert B\vert}2}\kern4cm\\
\kern2cm
\cdot
\bigg(\prod _{a\in A\cap B} ^{}x_a\bigg)^{2m+2}
\bigg(\prod _{a\in A^c\cap B^c} ^{}x_a^{-1}\bigg)^{2m+2}
\bigg(\prod _{a\in A\cap B^c} ^{}x_b\bigg)^{-1}
\bigg(\prod _{b\in A^c\cap B} ^{}x_b\bigg)\\
\cdot
V(A)V\left((A^c)^{-1}\right)R\!\left(A,(A^c)^{-1}\right)
V(B)V\left((B^c)^{-1}\right)R\!\left(B,(B^c)^{-1}\right).
\endmultline\tag\AV$$
We now fix disjoint subsets $A'$ and $B'$ of $[n]$ and extract the
coefficients of 
$$\bigg(\prod _{a\in A'} ^{}x_a\bigg)^{2m+n}
\bigg(\prod _{b\in B'} ^{}x_b^{-1}\bigg)^{2m+n}
$$
in (\AU) respectively in (\AV) (in the same sense as in
Section~4). If $\vert A'\vert+\vert B'\vert$ has
the same parity as $n$, then in (\AU) only the first sum contributes,
while in (\AV) it is only terms where the cardinality of the symmetric
difference $A\triangle B$ is even. One can then see in the same way as
Theorem~1 followed from Lemma~1, that the equality of corresponding
contributions follows from Lemma~2. On the other hand, if
$\vert A'\vert+\vert B'\vert$ has parity different from the parity of $n$, 
then in (\AU) only the second sum contributes,
while in (\AV) it is only terms where the cardinality of the symmetric
difference $A\triangle B$ is odd. In this case, the corresponding 
equality follows from Lemma~3.\quad \quad \qed
\enddemo

\proclaim{Theorem 4}For any non-negative integers $m$ and $n$, we have
$$\multline 
s_{((2m)^n)}(x_1,x_1^{-1},x_2,x_2^{-1},\dots,x_n,x_n^{-1})+
s_{((2m)^{n-1})}(x_1,x_1^{-1},x_2,x_2^{-1},\dots,x_n,x_n^{-1})\\
=
\sp_{(m^n)}(x_1,x_2,\dots,x_n)\,
o^{even}_{(m^n)}(x_1,x_2,\dots,x_n).
\endmultline
\tag\AW$$
\endproclaim

\demo{Sketch of proof} 
In a similar manner as we saw that the proof of Theorem~2 follows by
``replacing $m$ by $m+\frac {1} {2}$" in the proof of Theorem~1, the proof of
Theorem~4 follows by ``replacing $m$ by $m-\frac {1} {2}$" in the proof of
Theorem~3.\quad \quad \qed
\enddemo

Again, using (\ACb) and (\ACc), Theorems~3 and 4 allow for a uniform
statement, namely as
$$\multline 
\bigg(\prod _{i=1} ^{n}(x_i^{1/2}+x_i^{-1/2})\bigg)
\Big(s_{(M^n)}(x_1,x_1^{-1},x_2,x_2^{-1},\dots,x_n,x_n^{-1})\\
\kern4cm
+
s_{(M^{n-1})}(x_1,x_1^{-1},x_2,x_2^{-1},\dots,x_n,x_n^{-1})\Big)\\
=\so_{\big((\frac {M+1} {2})^n\big)}(x_1,x_2,\dots,x_n)\,
o^{even}_{\big((\frac {M} {2})^n\big)}(x_1,x_2,\dots,x_n).
\endmultline
\tag\AX$$

\medskip
In a similar vein as in Section~5, by specializing $x_i=1$,
$i=1,2,\dots,n$, in Theorems~3 and 4, we are able to derive potentially
interesting combinatorial interpretations. If we perform these
specializations in (\AT), then, by the combinatorial facts explained 
in Section~5, we obtain the identity
$$
\PP(2m+1,n,n) + \PP(2m+1,n-1,n+1) = \SPP(2m+2,n,n)\cdot\TCPP(2m,n,n),
\tag\AY$$
while the same specialization in (\AW) yields
$$
\PP(2m,n,n) + \PP(2m,n-1,n+1) =
\TCPP(2m,n+1,n+1)\cdot\SPP^*(2m,n,n),
\tag\AZ$$
where $\PP(A,B,C)$ denotes the number of plane partitions contained 
in the $A\times B\times C$ box (or, equivalently, the number of
rhombus tilings of a hexagon with side lengths $A,B,C,A,B,C$,
$\SPP(M,N,N)$ and $\TCPP(M,N,N)$ have the same meaning as in
Section~5, and where the $^*$ in $\SPP^*(2m,n,n)$ means that the
symmetric plane partitions in consideration satisfy the ``even 
central terraces
in the interior'' condition explained at the end of Section~5. 

Clearly again, the factorizations (\AY) or (\AZ) could be readily verified 
directly by using the product formulas for the combinatorial
quantities involved. However, it would have been difficult to see that
they exist without having first considered the factorization identities for
classical group characters given by Theorems~3 and 4, respectively.

\head \smc Acknowledgements \endhead

We are indebted to Ron King and Soichi Okada for extremely
insightful comments on the first version of this manuscript, which
helped to improve its contents considerably.

\enddocument